\input amstex   
\documentstyle{amsppt} 
\leftheadtext{Gerstenhaber, Giaquinto, Schaps}   
\rightheadtext{Donald--Flanigan problem for finite reflection
groups}
\magnification1200   
\parindent 20pt   
\NoBlackBoxes   

\topmatter   
\title The Donald--Flanigan problem for finite reflection groups  
   \endtitle 
\author Murray Gerstenhaber, Anthony Giaquinto, and Mary E.
Schaps 
      \endauthor 
\address Dept. of Mathematics, University of Pennsylvania, 
Philadelphia, PA 19104-6395 
      \endaddress 
\email mgersten\@math.upenn.edu     
 \endemail 
\address Dept. of Mathematics and Computer Science, Loyola
University of Chicago, Chicago, IL 60626--5311
     \endaddress
\email tonyg\@math.luc.edu
     \endemail
\address Dept. of Mathematics and Computer Science, Bar Ilan 
University, Ramat-Gan 52900, Israel 
      \endaddress 
\email mschaps\@bimacs.cs.biu.ac.il 
      \endemail 
\date August 4, 2000
     \enddate
\dedicatory To the memory of Mosh\'e Flato, z''l \enddedicatory
\abstract 
The Donald--Flanigan problem for a finite group $H$ and
coefficient ring $k$ asks for a deformation of the group algebra
$kH$ to a separable algebra. It is solved here for dihedral
groups and Weyl groups of types $B_n$ and $D_n$ (whose rational
group algebras are computed), leaving but six finite reflection
groups with solutions unknown. We determine the structure of a
wreath product of a group with a sum of central separable
algebras and show that if there is a solution for $H$ over $k$
which is a sum of central separable algebras and if $S_n$ is the
symmetric group then i) the problem is solvable also for the
wreath product $H\wr S_n = H \times\dots \times H\, \text{($n$
times)} \rtimes S_n$ and ii) given a morphism from a finite
abelian or dihedral group $G$ to $S_n$ it is solvable also for
$H\wr G$. The theorems suggested by the Donald--Flanigan
conjecture and subsequently proven follow, we also show, from a
geometric conjecture which although weaker for groups applies to
a broader class of algebras than group algebras. 
\endabstract
\endtopmatter  
\document

\subhead 1. Introduction: The Donald--Flanigan problem
\endsubhead
In this paper we solve the Donald--Flanigan problem for a large
class of groups including almost all finite reflection groups, 
and in the process give a simple construction of the rational
group rings (effectively, of all irreducible representations)
of the Weyl groups of types $B_n$ and $D_n$ (\S\S8-11). 

The Donald--Flanigan conjecture was one of the most intriguing of
those suggested by algebraic deformation theory because it sought
to relate the behavior of the group algebra of a finite group in
characteristic $p$ with that of its complex group algebra. If $G$
is a finite group and $k$ a commutative unital ring in which the
order of $G$ is invertible, then Maschke's theorem asserts that
the group algebra $kG$ is separable; in particular, $\Bbb CG$ is
a direct sum of matrix algebras. By contrast, if $k$ is a field
of characteristic $p$ dividing $\#G$ then $kG$ always has a  
non-trivial radical. Donald and Flanigan conjectured that for
this `modular' case the group algebra $kG$ can always be deformed
to a separable algebra, \cite{DF}. In the examples they
exhibited they observed that the separable algebra they
constructed was a direct sum of copies of the coefficient field
and thus resembled the complex group algebra. Schaps
suggested that one should look for a deformation of $kG$ with
matrix blocks in natural bijection with those of $\Bbb CG$,
corresponding blocks having the same dimensions. This concept is
made precise in the definition of a global solution.

The Donald--Flanigan conjecture fails for the 8-element
quaternion group
$\{\pm 1, \allowmathbreak \pm i, \allowmathbreak \pm j, \pm k\}$,
\cite{GGi}, but the problem remains of determining when
$kG$ can be deformed to a separable algebra of the ``right form''
(one resembling the complex group algebra). The result is
known for several classes of groups and algebras. (In the algebra
case, the problem is again to deform it to a separable algebra.)
These include abelian groups and symmetric groups $S_n$,
\cite{GSps}, (where in both cases there is a global solution
connecting the behavior at all primes), blocks with dihedral
defect group \cite{ES}, and blocks with abelian normal defect
group, \cite{MS}. Paradoxically, certain implications of the
Donald--Flanigan conjecture, cf. \cite{GG}, were subsequently
verified, in particular that in the modular case $kG$ always has
a non-inner derivation, \cite{FJL}. The paradox might be resolved
by the fact that this also follows (\S5) from a more geometric
conjecture which, although weaker for groups, is more general.

Finite Coxeter groups are identical with finite reflection
groups; we may use the terms interchangeably. There are four
infinite classes of irreducible ones (i.e., which are not direct
products), namely the Weyl groups of type $A_n (=
S_{n+1}), B_n, D_n$ and the dihedral groups. (For a general
reference and excellent exposition, cf. \cite{H}.) In addition,
there are four exceptional Weyl groups, $F_4, E_6, E_7, E_8$ and
two non-crystallographic groups, $H_3$ of order 120 (the symmetry
group of the icosahedron) and $H_4$ of order 14,400 (the symmetry
group of a regular 120-sided solid in $\Bbb R^4$). One can omit
$G_2$ because it is identical with the 12-element dihedral group
$\Cal D_6$.  Here we give another solution for the 
Donald--Flanigan problem for the dihedral groups (extending
certain work of Erdmann and Erdmann--Schaps, cf. \cite{ES}), and
using our earlier solution for $S_n$ we solve it for the groups
$B_n$ and $D_n$.  Thus there are only six finite Coxeter groups
for which the Donald--Flanigan problem remains open.  

An essential step in the solution is the determination of the
structure of a wreath product of a direct sum of central
separable algebras with a a finite group. Using this we show that
if $H$ is a group and $k$ a ring for which the Donald--Flanigan
problem has a solution which is a sum of central separable
algebras then i) the problem is solvable for the wreath product
$H\wr S_n = H \times\dots \times H \text{ ($n$ times) } \rtimes
S_n$ and ii) given a morphism from a finite abelian group or
dihedral group $G$ to $S_n$ it is solvable also for $H\wr G$.
While not important for us here, the condition on the solution
for $H$ will later be relaxed somewhat, since it is enough that
$kH$ be deformed to an algebra $A$ which becomes a direct sum of
central separable algebras after suitable extension of the
coefficient ring. For this it is sufficient that the center of
$kH$ be faithfully projective and of constant rank over the new
coefficient ring, an extension of $k$ resulting from the
deformation process.

Determining the structure of $H\wr S_n$ (and more generally, of
$H\wr G$ where $G$ operates by permutation of the factors of
$H^{\times n}$) parallels the classical problem of finding the
irreducible representations of a semidirect product $N\rtimes G$
in which the normal subgroup is abelian. The basic work
(fundamental to particle physics) is due to Wigner, who
computed the irreducible representations of the Poincar\'e
group (semidirect product of the Lorentz group and $\Bbb R^4$),
cf. \cite{St}, \S\S3.8, 3.9. Wreath products are special cases of
semidirect products but here the coefficient rings are generally
not fields and $H$ need not be commutative. Nevertheless, it may
be instructive to compare our treatment with that of Sternberg,
\cite{St}, particularly our example of \S3, where $H$ is
commutative, with that in \cite{St} (pp. 139--142) of the eight
element dihedral group $C_4 \rtimes C_2$, where $C_n$ is the 
$n$-element cyclic group. The Weyl group $D_n$ is only a
semi-direct product and not a wreath product, but the
Donald--Flanigan problem for it can be reduced to the wreath
product case. In the course of solving the Donald--Flanigan
problem for $B_n$ and $D_n$ we must, in effect, compute their
rational group algebras, which are given here explicitly.

Many related problems remain. In particular, suppose that a
finite group $G$ acts as automorphisms (possibly trivially) of a
finite group $H$, so that we can form the semidirect product
$H\rtimes G$. Then $G$ also operates as automorphisms of the
group algebra $kH$ for any coefficient ring $k$. An {\it
equivariant deformation} of $kH$ is one in which the operations
of $G$ continue to be automorphisms of the deformed algebra. 
When can an equivariant solution to the Donald--Flanigan problem
for $H$ be extended to one of the semidirect product?  More
generally, suppose that a finite group $G$ acts on a separable
$k$-algebra $A$ and let $A\#kG$, the {\it smash product}, be the
algebra with underlying $k$-module $A\otimes kG$ and
multiplication given by $(a\otimes\sigma) (b\otimes\tau) =
a\sigma(b)\sigma\tau;\, a,b\in A, \sigma,\tau\in G$. (This is a
special case of the usual smash product in Hopf theory, cf.
\cite{M}, Ch. 7, for note that $kG$ is a Hopf algebra with
$\Delta h = h\otimes h, h\in H$.) For example, 
$C_2\times \cdots \times C_2 \text{( $n$ times)} = C_2^{\times
n}$ is acted upon in a natural way not only by $S_n$ but by
$S_{n+1}$; the semidirect product is the Weyl group $D_n$. At
present we do not even know if under this operation $\Bbb
F_2C_2^{\times n}$ has an equivariant deformation to a separable
algebra. It seems unlikely. On the other hand, an equivariant
deformation may be more than one needs. It would be sufficient
that the action of the group deform simultaneously with the
structure of the algebra in such a way that the group continues
to act as automorphisms of the deformed algebra. This is what
happens in our solution to the Donald--Flanigan problem for $D_n$
but the behavior at the prime $p = 2$ is different from that at
other primes.

As preliminaries we reexamine what it means to deform an algebra
to a separable algebra and the concept of a global deformation
introduced in \cite{GSps}.

\subhead 2. Tempered successive deformations \endsubhead Recall
that a one-parameter family of deformations of an algebra $A$
over a ring $k$ is a $k[[t]]$-algebra $A_t$ whose underlying
module is $A[[t]]$ and which reduces modulo $t$ to the original
algebra $A$. It follows that multiplicaton in $A_t$ can be
written in the form $a*b = ab + tF_1(a,b) +t^2F_2(a,b) + \dots$,
where the $F_i$ are $k$-bilinear maps from $A$ to $A$ which
tacitly are extended to be $k[[t]]$-bilinear. Assume that 
 $A$ is free and of finite rank as a $k$-module.
Then as a $k[[t]]$-module $A_t$ is also free on the original
generators. One can similarly define multiparameter deformations
$A_{t,u,\dots,v}$ of $A$ in which the ``$*$'' multiplication is
given by power series in several variables. But no matter how
many parameters are used, if $A$ is not separable then no such
deformation can be separable, for $A$ remains a quotient of the
deformed algebra and a quotient of a separable algebra is always
again separable. 

For $k$ a field, saying that ``$A$ has been deformed to a
separable algebra'' tacitly means that there is some
multiplicatively closed subset $S$ of $k[[t,u,\dots,v]]$ such
that $S^{-1}A_{t,u,\dots,v}$ is separable over the new
coefficient ring $R = S^{-1}k[[t,u,\dots,v]]$. Here $R$ is free
over $k$ but when $k$ is only a (commutative, unital) ring this
seems too strong a condition; it should be sufficient (as in all
our examples) that $R$ be projectively faithful over $k$
although, of course, not of finite rank. This insures, in
particular, that $k$ remains a subring of $R$ and that a prime of
$k$ is invertible in $R$ if and only if it is already so in $k$.
Again, if $k$ is a field and there is only one parameter $t$ then
it is sufficient to invert $t$; the coefficient ring $R$ then
becomes the Laurent series field $k((t))$. Using the separability
idempotent, one can see that when $k$ is arbitrary and there are
several parameters it is still necessary to invert only a single
element of $k[[t,u,\dots,v]]$. For any single element of
$S^{-1}A_{t,u,\dots, v}$ is already contained in the algebra
obtained by a single such inversion and the separability
idempotent only involves a finite number of elements. But even in
the case of a single parameter, it may not be simply $t$ that one
wants to invert but some polynomial in $t$ (whose constant term
is not a unit), cf. \cite{GSps}. 

Suppose now that we try to deform an algebra $A$ to a separable
algebra in two stages: first reducing the inseparability by
deforming to some $A_t$ and forming $B=S^{-1}A_t$ where $S$ is
generated by some single element of $k[[t]]$, and then by
deforming $B$ to some $B_u$ and inverting an element of the new
cofficient ring. Can the same be done by first performing some
two-parameter deformation $A_{t,u}$ of $A$ and then inverting an
element of $k[[t,u]]$? We need this in order to justify deforming
an algebra to a separable one in successive steps. To understand
the reason for caution, suppose that $A$ is a finite-dimensional
algebra over a field $k$, that we deform it, and then specialize
the deformation parameter. If $\dim_k A = n$ and we have chosen a
basis $a_1,\dots,a_n$ then the algebra can be described by its
multiplication constants relative to this basis: $a_ia_j = \sum
c_{ijl}a_l$. The constants $c_{ijl}$ may be viewed as determining
a $k$-point in the variety $\operatorname{alg}_n$ of structure
constants of $n$-dimensional associative algebras. The deformed
algebra has structure constants $c_{ijl}(t)$ lying in $k[[t]]$
which reduce at
$t = 0$ to the original constants. These in general will be
transcendental over $k$ and define a subvariety in the same
component of $\operatorname{alg}_n$ as the original. The result
of specializing and then performing another deformation can,
however yield a point in a different component, since the
specialization may lie on the intersection of two components. But
intuitively, if one deforms and then, starting from a
generic point of the original deformation deforms again, one must
remain on the same component. 

So suppose that $A$ is free as a $k$-module with basis
$a_1,\dots,a_n$, that we have deformed it to $A_t$ and have
formed $B=S^{-1}A_t$ where $S$ is generated by a single element
$f(t) \in k[[t]]$. Let $B$ now be deformed to $B_u$.
Set $S^{-1}k[[t]] = R$. We may suppose that we have a basis
$b_1,\dots,b_n$ of $B$ with structure constants given by
$b_{\lambda}b_{\mu} = \sum \gamma_{\lambda\mu\nu}b_{\nu}$ where
the  $\gamma_{\lambda\mu\nu}$ initially lie in $R$ but after
deformation are elements of $R[[u]]$. The problem is that the
coefficients of these power series in $u$ may contain negative
powers of $f(t)$ and the negative powers may be unbounded. The
hypothesis we make about the deformation $B_u$ is that there is
some fixed $N$ such that all $f(t)^N\gamma_{\lambda\mu\nu}(u)$
lie in $R[[t,u]]$. Now we can write each of the $b_{\lambda}$ as
a linear combination of the $a_i$ with coefficients in $R$.
Again, these coeficients will generally involve negative powers
of $f(t)$. As there are only finitely many coefficients, the
hypothesis insures that we can write out the multiplication in
$B_u$ in the form $a_i*a_j = \sum c_{ijl}(t,u)a_l$ where there is
some fixed $N'$ such that all $f(t)^{N'}c(t,u) \in k[[t,u]]$. But
all of the $c(t,0)$ lie in $k[[t]]$. Therefore, replacing $u$ by
$f(t)^{N'}v$ we have an ordinary two-parameter deformation
$A_{t,u}$ of $A$ over $k[[t,u]]$. A second deformation with the
foregoing boundedness properties will be called {\it
tempered} (relative to the first). If, finally, inverting an
element of $R[[u]]$ with bounded powers of $f(t)$ in the
denominators of the coefficients makes $B_u$ separable, then
there is an element $F(t,u) \in k[[t,u]]$ inverting which makes
$A_{t,u}$ separable. To return to a one-parameter family of
deformations we can now replace $u$ by some element of $k[[t]]$
such that $F(t,u) \ne 0$.

\subhead 3. An example\endsubhead
As a simple illustration, let $C_2 = \{1, a, a^2=1\}$ be the two
element group and $G = C_2\wr C_2 = (C_2 \times C_2)\rtimes C_2$.
The right $C_2$  operates by interchange of the two left factors.
Denoting its elements by $\{1,\sigma\}$, we will write the 8
elements of the group as $\{(1,1), (a,1), (1,a), (a,a),
(1,1)\sigma, (a,1)\sigma, (1,a)\sigma, (a,a)\sigma\}$ where, for
example, $(a,1)\sigma\cdot (1,a)\sigma = (a,1)(a,1)\sigma^2 =
(1,1)$, the unit element of $G$. Now letting $k = \Bbb F_2$, we
wish to deform the group algebra $\Bbb F_2G$ to a separable
algebra. We can deform $\Bbb F_2C_2$ to a separable algebra by
setting $a^2 = ta + 1+t$; then $(1+a)^2 = t(1+a)$ so after
inverting $t$ we have orthogonal idempotents $e = (1+a)/t$ and $f
= 1+e$. Doing the same to both factors in $C_2 \times C_2$ (and
noting that if $G_1, G_2$ are finite groups then $k(G_1\times
G_2) = kG_1 \otimes kG_2$) we can perform a first deformation of
$\Bbb F_2G$ to an algebra generated by the four orthogonal
idempotents $e\otimes e, f\otimes e, e\otimes f, f\otimes f$ and
a ``switch'' element $\sigma$ with $\sigma^2 = 1$ such that
$\sigma(x\otimes y)\sigma = (y\otimes x)$ for all $x\otimes y \in
\Bbb F_2C_2 \otimes \Bbb F_2C_2$. The resulting algebra is still
in a natural sense a wreath product. Let $A$ be the algebra to
which $\Bbb F_2C_2$ has been deformed. Its coefficients are now
in $\Bbb F_2(t)$. We still have $C_2$ operating on $A\otimes A$
by interchange of the tensor factors. Denote the algebra
resulting from this first deformation by $A\wr C_2$. Further, $A$
is separable over $\Bbb F_2(t)$ (part of the original
inseparability has been removed) and is a direct sum of two
subalgebras, $A = Ae\oplus Af$, each of which is trivially
central separable. 
While $A\otimes A$ has four central primitive idempotents, $A\wr
C_2$ is a direct sum of 3 subalgebras corresonding to the orbits
of these idempotents under $C_2$, namely $\{e\otimes e\},
\{f\otimes f\}$, and $\{e\otimes f, f\otimes e\}$. The isotropy
group of the last orbit is reduced to the identity element of
$C_2$ and the orbit gives rise to a four dimensional summand of
$A\wr C_2$ spanned over $\Bbb F_2(t)$ by $\{e\otimes f, f\otimes
e, (e\otimes f)\sigma, (f\otimes e)\sigma\}$. It is easy to check
that this is isomorphic to the $2\times 2$ matrix algebra
$M_2(\Bbb F_2(t))$ and hence is central separable over $\Bbb
F_2(t)$. The isotropy groups of the other orbits are non-trivial,
being
in fact all of $C_2$. The subalgebra corresponding to $\{e\otimes
e\}$ is spanned by $\{e\otimes e, (e\otimes e)\sigma\}$ and
should be viewed as isomorphic to $M_1(\Bbb F_2(t))\otimes\Bbb
F_2(t)C_2$, where $C_2$ here is the isotropy group of the orbit,
and similarly for the orbit of $\{f\otimes f\}$. This will be
generalized in \S\S8, 9. We can now perform a second
deformation, treating the direct summands individually, so that
the whole algebra now becomes separable over the new coefficient
ring. The non-trivial matrix summand is already central separable
so it will be left unchanged except for the necessary extension
of the coefficient ring. In this second deformation we can not
simply set $\sigma^2 = u\sigma + 1+u$ (since there is to be no
change in the matrix part); in the first summand we set
$[(e\otimes e)\sigma]^2 = u(e\otimes e)\sigma + (1+u)(e\otimes
e)$ and similarly with $f$ in the second. 

It is easy to see that this second deformation is tempered
relative to the first. In fact, $\Bbb F_2(t)C_2$ is actually
defined over $\Bbb F_2$; it is obtained from $\Bbb F_2C_2$ just
by extension of coefficients. So the only way that $t$ enters
into the second deformation is in the choice of the basis
for the algebra obtained by the first deformation. Since there
are only a finite number of basis elements, the condition of
being tempered will automatically be satisfied whenever they all
have coefficients (relative to the original basis) in some
$S^{-1}k[t]$ where $S$ is generated by a single $f(t) \in k[t]$.
Nevertheless, it may be useful in this example to write
explicitly the final result as a two parameter deformation over
$\Bbb F_2[[t,u]]$ (in fact, over $k[t,u]$) followed by inversion
of $tu$. Since now $\sigma*\sigma$ is no longer the unit element
of the algebra, we compute it explicitly. We can write $\sigma =
[(e+f)\otimes(e+f)]\sigma = [e\otimes f + f\otimes e]\sigma
\oplus (e\otimes e)\sigma \oplus (f\otimes f)\sigma$ relative to
the decomposition after the first deformation. The second
deformation has respected this decomposition, so we can square
each direct summand separately. The square of the first is
$e\otimes f + f\otimes e$, that of the second is $(1+u)(e\otimes
e) + u(e\otimes e)\sigma$, and similarly for the third with $f$
replacing $e$. Thus $\sigma*\sigma = (e+f)\otimes(e+f)
+u(e\otimes e + f\otimes f)(1+\sigma)$. Writing out
$\sigma*\sigma$ in terms of the original basis elements the first
summand remains unchanged (being just the unit element of both
the original and the deformed algebra) but the second summand
must be written as $t^{-1}u(a\otimes 1 + 1\otimes a +
t\cdot1\otimes 1)(1+\sigma)$. So we need here to replace $u$ by
$tv$, which in fact will work for all the other products. With
this we have a true two-parameter  deformation of the original
algebra (the parameters now being $t$ and $v$) which becomes
separable after inversion of $t$ and of $u =tv$, or simply after
inversion of $tv$. To get a  one-parameter deformation to a
separable algebra we could now set $u=t^2$. (Later we will
actually use the deformation given by $a^2 = (q-q^{-1})a + 1$
with $q = 1+t$; the results are the same.)

While we must repeatedly use successive deformations, in the
cases we consider it will be evident, as it is here, that the
second deformation is tempered, so we may simply omit the
discussion of that fact. 

\subhead 4. Global solutions and split global solutions
\endsubhead In the modular case we would like to deform $kG$ to a
separable algebra which is a direct sum of matrix algebras whose
summands are in one-one correspondence with the matrix blocks of
the complex group algebra $\Bbb CG$, with corresponding blocks
having the same dimension. More precisely, let $k$ now be a
field.  Write $kG = A$, denote the deformed algebra by $A_t$, and
supppose that it becomes separable when coefficients are extended
to $k((t))$. It need not be a sum of matrix algebras but will
become one when coefficients are extended to the algebraic
closure of $k((t))$, i.e., when one forms $A_t
\otimes_{k[[t]]}\overline{k((t))}$. But in general there need be
no relation between the matrix blocks of this algebra and those
of $\Bbb CG$. For example, $G$ may be abelian but $kG$ may have
non-commutative separable deformations. Suppose, however, that 
we have a {\it global solution} to the Donald--Flanigan problem 
for $G$ in the sense of \cite{GSps}: a deformation $A_t$ of the
integral group ring $\Bbb ZG$, together with a multiplicatively
closed subset $S$ of $\Bbb Z[[t]]$ which does not meet the ideal
generated by any rational prime dividing $\#G$, such that
$S^{-1}A_t$ is separable over $S^{-1}\Bbb  Z[[t]]$. If such a
global solution exists then one may assume that $S$ consists of
the powers of a single element $s$. More important, one can
reduce $S^{-1}A_t$ modulo any rational prime $p$ not dividing
$\#G$. The image $\bar s$ of $s$ must be of the form
$t^m\epsilon$, where $\epsilon$ is a unit of $\Bbb 
F_p[[t]]$. Since a quotient of a separable algebra is separable,
the result will be a deformation of $\Bbb F_pG$ which becomes
separable when coefficients are extended to $\Bbb F_p((t))$. 
Denote the resulting separable algebra by $A_t(p)$. Now if
$S^{-1}A_t$ is already split, i.e., a direct sum of matrix
algebras, then not only must the same be true of each $A_t(p)$ 
but the correspondence between the matrix blocks of $S^{-1}A_t$
and those of $A_t(p)$ is simply that those of the latter are the
quotients of those of the former. Corresponding blocks then
certainly also have the same dimensions. Since one can embed
$\Bbb Z[[t]]$ in $\Bbb C$ it also becomes clear that the blocks
are the same as those of $\Bbb CG$. If $S$ does contain some
primes $p_1, p_2, \dots$ dividing $\#G$ then we say that we have
a {\it global solution away from} $p_1, p_2, \dots$. An arbitrary
deformation of $kG$ to a separable algebra will be called a {\it
weak solution} to the Donald--Flanigan problem. Different weak
solutions are sometimes possible. For example $\Bbb F_2(C_2
\times C_2)$ is deformable both to a direct sum of four copies of
the new coefficient ring and to a $2 \times 2$ matrix ring over
the coefficient ring. 
A principal result of \cite{GSps} is that there does exist a
split global solution to the Donald--Flanigan problem for the
symmetric groups $S_n$. In general, however, a global solution if
it exists will not be split. For letting $t \to 0$, the existence
of a split global solution would imply that all irreducible
representations of $G$ are rational, something true of $S_n$ but
not in general. This problem would disappear if we had a positive
answer to the following question. If $A$ is a separable algebra
over a domain $R$, is there always a finite integral extension
$\hat R$ of $R$ over which $A$ is split? For if $\Bbb ZG$
is deformed to a separable algebra over some $R =
S^{-1}Z[[t,u,\dots,v]]$, where $S$ is a multiplicatively closed
set not intersecting the ideal generated by any rational prime
dividing $\#G$, then the same will still be true if the
coefficients are extended to $\hat R$; if this splits the
separable deformation then we can proceed as before. Known
results  do not get us quite this far; see \S10. Note that
for our purposes one can assume that $R$ has characteristic zero
and even that all rational primes not dividing $\#G$ are already
invertible in $R$, but it does not seem that these assumptions
should be necessary. The original definition of a global solution
is also too restrictive in that it requires that we start with
the integral group ring $\Bbb ZG$ of the group $G$. It is useful
to broaden the definition by allowing replacement of $\Bbb Z$ by
some finite integral extension $\Cal O$. This need not be the
full ring of integers of some number field but we may generally
assume that it is.  When $G$ and a single $p|\#G$ are given then
a global solution away from all other primes dividing $\#G$ will
be called a {\it local solution} at $p$. If either it is already
split or the conjecture above holds then this gives a canonical
correspondence between the matrix blocks of $\Bbb CG$ and those
of the global solution. Suppose now that we have local solutions
at several primes dividing $\#G$. Then their blocks must be in
natural correspondence with each other, corresponding blocks
having the same dimensions, since for each prime they are in
correspondence with those of $\Bbb CG$. The requirement that the
local or global solution be split can be eased. All one needs is
that it be a direct sum of Azumaya, i.e., central separable
algebras. Such an algebra, if it has constant
rank at each prime ideal of its coefficient ring is a just a
twisted form of a matrix algebra, cf, e.g., \cite{KO}. We return
to this in \S10.

The global solution (with $\Cal O$ just $\Bbb Z$ itself) for the
Donald--Flanigan problem for $S_n$ given in \cite{GSps} is
essentially its Hecke algebra, $H_n(q)$. Setting $\Bbb Z_q = \Bbb
Z[q, q^{-1}]$, this is the free module over $\Bbb Z_q$ with basis
elements $T_w$ indexed by the elements  $w \in S_n$  and
multiplication given as follows: The {\it length} $\ell(w)$ is
the number of factors in a shortest expression of  $w$ as a
product of generators $s_i := (i,i+1), i = 1,\dots,n-1$ of
$S_n$.  Now set (i) $T_sT_w = T_{sw}$ if $s = (i,i+1)$ for some
$i$ and  $w\in S_n$ is an element with $\ell(sw) > \ell(w)$, and
(ii) $T_s^2 = (q-q^{-1})T_s+1$.  This implies that $T_sT_w =  (q-
q^{-1})T_w + T_{sw}$ when $\ell(sw) <\ell(w)$. (The same
definition extends to any Coxeter group, and in particular to any
finite reflection group where instead of the transpositions $s_i
=(i,i+1)$ one takes the basic generators $s$.)  Writing $1+t$ for
$q$  one sees that $H_n(1+t)$ is a deformation of $\Bbb ZS_n.$
Set $i_q:= (1-q^i)/(1-q)$ and similarly 
$i_{q^2}:= (1-q^{2i})/(1-q^2)$; set $n_{q^2}! :=
n_{q^2}(n-1)_{q^2}\dots 2_{q^2}$ and $\Bbb  Z_{q,n}: = \Bbb Z[q,
q^{-1}, 1/n_{q^2}!]$.  A main result of \cite{GSps} was that over
$\Bbb Z_{q,n}$ the Hecke algebra  $H_n(q)$ becomes a direct sum
of matrix algebras. It follows that $H_n(1+t)$ together with the
multiplicatively closed subset of $\Bbb Z[[t]]$ generated by
$n_{q^2}! = n_{(1+t)^2}!$ is a global solution to the 
Donald--Flanigan problem for the symmetric group. Setting $q=1$
or equivalently $t=0$ yields the special case that all the
irreducible complex representations of the symmetric group are
actually defined over $\Bbb Q$. This is the reason that no
extension of $\Bbb Z$ was needed for $G = S_n$. In general
$i_{q^2}j_{q^2} \ne (ij)_{q^2}$, so adjoining the inverse of
$m_{q^2}$ generally does not bring with it the adjunction of
$i_{q^2}$ for any factor $i$ of $m$. In particular, the
multiplicatively closed set $S$ above is not generated by
$n!_{q^2}$ ($\ne n_{q^2}!$). However, if $m$ is even, say $m =
2r$, then $(2r)_{q^2} = (1+q^r)r_{q^2}$, so inverting
$(2r)_{q^2}$ also inverts $r_{q^2}$. To make the Hecke algebra of
$S_n$ separable we therefore essentially adjoined to its original
ring of definition $\Bbb Z_q$ the inverses of all $i_{q^2}$ for
$i = 2,\dots,n$. There is no proof that all of these adjunctions
are necessary, but we suspect that they are. The numbers $2,
\dots, n$ are are also the ``degrees'' of $S_n$, i.e., degrees of
its basic invariant polynomials. This suggests the following
refinement of a conjecture in \cite{GSps}. Let $W$ be a finite
Coxeter group and $H_W(q)$ be its Hecke algebra with coefficient
ring $\Bbb Z_q = \Bbb Z[q, q^{-1}]$. If $S$ is the
multiplicatively closed subset generated by all $i_{q^2}$ where
$i$ runs through the degrees of $W$ then $S^{-1}H_W(q)$ should be
separable. We do not know if the deformations of $\Bbb Z B_n$ and
$\Bbb Z D_n$ constructed here are in fact their Hecke algebras,
but they do have the following property.  Setting $n_{q^2}! =
n_{q^2}(n-1)_{q^2}\cdots 3_{q^2}2_{q^2}$ and $\Bbb Z_{q,n} = \Bbb
Z_q[1/n_{q^2}!]$ the deformed algebras become direct sums of
matrix algebras over $\Bbb Z_{q,n}$. Setting $q=1$ recaptures the
result that all the complex irreducible representations of $B_n$
and $D_n$ are rational. Since this is not the case for the
dihedral groups (whose degrees are $2,m$), even if the
conjecture is true then the resulting solution is not split. On
the other hand, since the product of the degrees is equal to the
order of $W$, it is consistent with the fact that the group
algebra $kW$ is separable over any ring $k$ in which $\#W$ is
invertible. If this conjecture is true then in principle one
could prove it by exhibiting the separability idempotent, but as
noted in \cite{GSps} there is so far no good formula for that
even in the known case where $W = S_n$. It may be the case more
generally that rings of invariant polynomials associated with a
finite group have some connection with the solvability of the
Donald--Flanigan problem for that group. The conjecture does not
assert that the Hecke algebra is the only global solution to the
Donald--Flanigan problem for a finite reflection group (although
it may in some sense be the best). For the Weyl groups of type
$B_n$ the degrees are $2, 4,\dots,2n$, but we shall see that one
can obtain a global solution as soon as one has one for $S_n$,
and that requires inverting only all $i_{q^2}$ for $i=1,\dots,n$.
However, as remarked above, if a ring contains the inverses of
all $(2i)_{q^2}$ for $i=1,\dots,n$ then it already contains the
inverses of all $i_{q^2}$.

Infinite Coxeter groups also have Hecke algebras which 
can be viewed as deformations of their integral group rings. If
coefficients are extended to the direct limit $\Bbb Z_{q,\infty}$
of the
$\Bbb Z_{q,n}$ then these algebras should possess some properties
similar to
separability. While $\Bbb Z_{q,\infty}$ is not a field, all the
``quantum integers'' $n_{q^2}$ have become invertible. While
almost a quantum version of the rationals, it is still possible
to reduce $\Bbb Z_{q,\infty}$ modulo any rational prime. Although
reduction modulo $t$ does give the rationals, it is certainly not
a deformation of the rationals

\subhead 5. Geometric rigidity and a generalized global problem
\endsubhead 
For this and the next section the coefficient ring $k$  will be a
field. We review some basic facts about jump deformations and
approximate automorphisms. 

Recall that a {\it jump deformation} $A_t$ of a $k$-algebra $A$
is one which is non-trivial and remains constant for generic $t
\ne 0$. More precisely, if $u$ is a second variable and
coefficients are extended to $k((t))[[u]]$ then there is an
isomorphism $A_t \cong A_{(1+u)t}$ which reduces to the identity
when $u\to 0$. That is, the trivial deformation of the 
$k((t))$-algebra $A_t$ is equivalent to the deformation
$A_{(1+u)t}$. (For the essential properties of jump deformations,
cf. \cite{G2,3; GSck, \S7}.) A jump deformation simultaneously
``breaks some symmetry'' of the algebra and ``destroys its own
infinitesimal''. The meaning of the latter statement is this.
Suppose that we have a deformation $A_t$ of $A$. An $n$-cocycle
$\zeta \in Z^n(A,A)$ can be {\it lifted} to an $n$-cocycle of
$A_t$ if there are cochains $\zeta_i \in C^n(A,A)$ such that
$\zeta_t = \zeta + t\zeta_1 +t^2\zeta_2 + \dots \in Z^n(A_t,
A_t)$. If $\zeta
+t\zeta_2 +\dots +\zeta_m$ is a cocycle modulo $t^{m+1}$ then
there is an obstruction cocycle $\eta \in Z^{n+1}(A,A)$ and we
must have $\eta = \delta\zeta_{m+1}$ for some $\zeta_{m+1} \in
C^n(A,A)$ in order for the construction to continue. Any cocycle
which is liftable to a coboundary is called a {\it jump} cocycle.
This includes, in particular all coboundaries in $A$ but there
may be non-trivial ones. One has $H^n(A_t,A_t)$ = (liftable  
$n$-cocycles)/(jump $n$-cocycles).  The infinitesimal of a jump
deformation is always a jump cocycle, so $\dim H^2(A_t,A_t) <
\dim H^2(A,A)$. It follows that $A$ can only undergo finitely
many jump deformations. 

For $n\ge 3$ every jump cocycle in $Z^n(A,A)$ is the obstruction
at some stage to lifting some $\zeta \in Z^{n-1}(A,A)$. In
characteristic zero this is true also for $n=2$. In this case,
therefore, if $F$ is the infinitesimal of a jump deformation then
$F$ is itself the obstruction to lifting some derivation $\phi
\in Z^1(A,A)$ to a derivation of $A_t$. Therefore $e^{t\phi}$ is
a formal one-parameter family of automorphisms of $A$ which can
not be lifted to $A_t$, so this symmetry is broken. In
characteristic $p > 0$ it is possible that after a jump
deformation every $\phi \in Z^1(A, A)$ remains liftable, but
there is still a symmetry that is broken. For example, in $A =
\Bbb F_2[x]/x^4$ a derivation is uniquely determined by its value
on $x$, which can be arbitrary. This is still the case if we
deform $A$ to $A_t = \Bbb F_2[x]/(x^4+tx^2)$ (a jump
deformation), so the dimension of the space of derivations has
not changed. We have, however, lost some ``approximate
automorphisms'' in the following sense. When we have a derivation
$\phi$ of an algebra $A$ of characteristic $p > 0$ it may not be
always be possible (as it was in characteristic zero) to
construct a full one-parameter family of automorphisms of the 
form $\Phi_t = \phi + t\phi_1 + t^2\phi_2 + \dots$. The largest
$m$ for which we can construct an automorphism of $A[t]/t^m$ of
the form $\Phi_t = \phi + t\phi_1 + \dots +t^{m-1}\phi_{m-1}$
with a fixed $\phi$ always has $m=p^r$ for some $s$. We call this
$m$ the order of the {\it approximate automorphism} $\Phi_t$ and
also the order of $\phi$. A jump deformation of $A$ reduces the
order of some derivation. In the example above, the order of the
derivation $\phi$ of $\Bbb F_2[x]/x^4$ sending $x$ to $1$ was
initially 4 but was reduced to 2 after the jump deformation. With
a jump deformation the order of some derivation is always
reduced; it can not be lifted to an approximate autmorphism of
the same order, and the infinitesimal of the deformation is the
obstruction to continuing the now truncated approximate
automorphism to higher order. In this sense, a jump deformation
always breaks some symmetry. If $H^1(A,A) = 0$ there can be no
jump deformations.

\remark{Remark} Following \cite{G2}, denote by
$\operatorname{Aut}_{m-1}A$ the group of automorphisms of
$A[t]/t^m$ of the form $\operatorname{id}_A +t\phi_1+ \dots
t^{m-1}\phi_{m-1}$
where the $\phi_i$
are 1-cochains of $A$, i.e., linear maps $A \to  A$ (which are
tacitly extended to be $k[t]$-linear). There is a canonical
monomorphism $\operatorname{Aut}_{m-1}A \to
\operatorname{Aut}_{mn-1}$ defined by replacing $t$ by $t^n$ and
considering the resulting polynomial as one of degree $mn-1$. The
direct limit $\varinjlim\operatorname{Aut}_mA$ contains a
subgroup $P$ consisting of the images of all elements in all
$\operatorname{Aut}_mA$ which have a prolongation to an element
of $\operatorname{Aut}_{m+1}A$. This subgroup is normal by the
basic theorem on the additivity of obstructions \cite{G2, Theorem
1}, and the quotient $\operatorname{Aut}_*A$ is abelian. It
consists of classes of obstructed approximate
automorphisms, since all full one-parameter families of
automorphisms beginning with the identity lie in $P$. This group
has a natural filtration, as does the space $RH^2(A,A)$ of
``restricted'' elements of $H^2(A,A)$ (=those which are the
obstructions to approximate automorphisms) and the associated
graded groups are isomorphic. The group of obstructed approximate
automorphisms is thus in a natural way finite dimensional. In
characteristic zero the Euler-Poincar\'e characteristic of an
algebra $A$, if it has one, is just $\sum (-1)^n\dim H^n(A,A)$.
If it does, then so does any deformation of $A$ and it is
invariant under deformation. In characteristic $p>0$ for this to
be true we shall probably have to replace $\dim H^1(A,A)$ with
the dimension of the space of approximate automorphisms.
\endremark 

The original weak Donald--Flanigan conjecture implied, by way of
the following theorem (\cite{GGr}) that if $p|\#G$ then there
exists an element $g\in G$ whose centralizer $C_G(g)$ contains a
normal subgroup of index $p$. 
\proclaim{Theorem 1} Let $A$ be a $k$-algebra which is not itself
rigid but which can be deformed to a rigid algebra. Then
$H^1(A,A)\ne 0$, i.e., $A$ has a non-inner derivation.
\endproclaim \demo{Proof} The infinitesimal of any deformation of
$A$ to a rigid algebra must be a jump 2-cocycle, since by
definition no further deformation of a rigid algebra is possible.
Since a jump 2-cocycle exists, $H^1(A,A) \ne 0$. $\qed$ \enddemo

Any deformation of a non-separable algebra $A$ to a separable one
is necessarily a jump deformation since separable algebras have
trivial cohomology and therefore are rigid. It follows that
$H^1(A,A)\ne 0$. However, for $A=\Bbb F_pG$ we have (cf. \cite{B,
Theorem 2.11.12}) an isomorphism of additive groups 
$$H^n(\Bbb F_pG, \Bbb  F_pG) \cong \bigoplus
H^n(C_G(g), \Bbb F_p)$$ 
where the operation on $\Bbb F_p$ is trivial and the sum is over
a set of representatives $g$ of conjugacy classes in $G$. So the
Donald--Flanigan conjecture implied that for some $g\in G$ one
has $H^1(C_G(g), \Bbb F_p) \ne 0$. But a derivation into a
trivial module is just a morphism, giving what was asserted.
Guided by this Fleischmann, Janiszczak and Lempken proved a
stronger  result \cite{FJL}: If $p|\#G$ then
there exists a $g \in G$ whose order is divisible by $p$ and
whose centralizer $C = C_G(g)$ has the property that its
commutator subgroup $C'$ does not contain the $p$-part of $g$.
Their proof reduces to the case where $G$ is simple and uses the
classification theorem for finite simple groups. (Publication of
\cite{FJL} preceded that of \cite{GGr} because of the greater
lead time for the latter.)

To resolve the paradox that the Donald--Flanigan conjecture fails
while its corollary in \cite{GGr} holds we propose a conjecture
which for finite groups is weaker than the Donald--Flanagan
conjecture but applies more generally and still implies the
statement in \cite{GGr}. Recall first the various concepts of
rigidity for an algebra $A$, cf. \cite{GSck}. The first, usually
called simply rigidity but more precisely {\it analytic rigidity}
says that every formal deformation of $A$ is equivalent to the
trivial deformation. This will certainly hold if $H^2(A,A) = 0$,
often called {\it absolute rigidity}. For the second, suppose for
the moment that $A$ has dimension $n < \infty$ over some
algebraically closed field $k$ and let $\Cal V =
\operatorname{alg}_n(k)\subset k^{n^3}$ be the variety of
structure constants of $n$-dimensional $k$-algebras. Now
$GL(n,k)$ operates on $\Cal V$ and hence on each of its
components, with orbits corresponding to isomorphism classes of
$n$-dimensional $k$-algebras. One calls $A$ {\it geometrically
rigid} if the corresponding orbit is a Zariski open set in the
component $\Cal V_A$ of $\Cal V$ containing it. In the
terminology of Wigner, every algebra in the (necessarily unique)
component of $A$ is a contraction of $A$. Equivalently, every
algebra which can be deformed to $A$ is a contraction of $A$;
this formulation no longer requires finite dimensionality. 
Analytic rigidity implies geometric rigidity and in
characteristic zero they are equivalent \cite{G3, GSck \S9}, but
not in characteristic $p>0$. In that case there can exist 
non-trivial formal deformations of $A$ with multiplication $a*b =
ab + tF_1(a, b) + t^2F_2(a,ab) + \dots$ in which $F_1$ is not a
coboundary but where the deformation becomes trivial when $t$ is
replaced by $t^m$ for some $m$; this $m$ is then necessarily a
power of $p$, cf. \cite{G3, GSck}. Such {\it restricted}
deformations can occur only when $A$ has some non-inner
derivations $\phi$, for the cocycle $F$ must be the obstruction
at some stage to constructing a formal automorphism of the form
$\Phi = \operatorname{id}_A + t\phi + t^2\phi_2 + \dots$ whose
infinitesimal is $\phi$.  The obstruction to a derivation is
always the infinitesimal of a restricted deformation. (See also
Skryabin, \cite{Sk}).

Suppose that an algebra $A$ is defined and is a free module of
finite rank over an integral extension $\Cal O$ of $\Bbb Z$ (so
the structure constants with respect to a suitable choice of
basis lie in $\Cal O$), and suppose that $A$ becomes rigid over
$\Bbb C$. In place of Donald--Flanigan we conjecture now that for
any prime $p$, $A \otimes_{\Cal O} \Bbb F_p$ can be deformed to a
geometrically rigid algebra. The conjecture relates the structure
of the component $\Cal V_A$ of $A$ over $\Bbb C$ to the structure
of the variety of structure constants in characteristic $p$ at
any point representing the reduction of $A$ mod $p$. If one has a
finite group $G$ then the hypothesis of having a free module  is
certainly satisfied (with $\Cal O$ just $\Bbb Z$ itself) for  the
integral group algebra $\Bbb ZG$. The conclusion is weaker  than
that of the Donald--Flanigan conjecture, for if $k$ is a field of
characteristic $p$ dividing $\#G$ then it says only that $kG$ can
be deformed to a geometrically rigid algebra, not a separable one
(all of whose cohomology in positive dimensions vanishes). This,
however, is enough to imply that $\Bbb F_pG$ has a non-inner
derivation. For if the geometrically rigid algebra has a
restricted deformation then we are done; otherwise it is
analytically rigid, the deformation is just as before a jump
deformation whose infinitesimal is the obstruction to some
derivation, so in any case the derivation exists. One is then led
again to the group-theoretic statement in \cite{GGr}. 

For group algebras over a field, the conjecture of the next
section would imply that they can actually be deformed to
geometrically rigid algebras, but to state it we must introduce a
broader concept of equivalence of deformations $*$ and $*'$ of a
$k$-algebra $A$. We will call them {\it effecively equivalent} if
there are integers $m$ and $n$ such that replacing $t$ by $t^m$
in the first and by $t^n$ in the second, the resulting
deformations become equivalent in the original sense of
\cite{G1}. For finite dimensional algebras $A$ over an
algebraicially closed field $k$ this can happen in a non-trivial
way only in characteristic $p > 0$, where its significance is the
following. Choosing a basis for $A$, one has new multiplication
constants given by $*$ and by $*'$ which in both cases are now
power series in $t$ reducing to the original multiplicaton
constants when $t=0$. Each set of constants may now be viewed as
a generic point of some subvariety of $k^{n^3}$; the deformations
are effectively equivalent precisely when these subvarieties
coincide. A deformation which is effectively equivalent to a
trivial deformation is one previously called restricted but a
better name may be {\it effectively trivial}. An {\it effective
jump} deformation will be one which is effectively equivalent to
a (non-trivial) jump deformation.

\subhead 6. Jump algebras \endsubhead 
There is an important class of algebras $A$ where (using known
deep results) it is easy to show that $H^1(A,A) \ne 0$. These are
those algebras of finite representation type or
representation finite algebras which are not already rigid.
(Clearly one must exclude, e.g., matrix algebras.) The reason is
that the only deformations which these admit are effectively
equivalent to jump deformations (and hence, in characteristic
zero, they admit only true jump deformations). In this section
the coefficient ring $k$ is still a field (although some
statements will obviously hold more generally), algebras will be
finite dimensional, and a deformation of a $k$-algebra $A$ will
be viewed as an algebra over the Laurent series field $k((t))$. 
A class of algebras is called {\it open} if every deformation of
an algebra in the class is again in the class. For algebras of a
fixed dimension $n$, an open class may be represented as a
Zariski open subset of the variety of structure constants of 
$n$-dimensional algebras. A theorem of Gabriel asserts that the
class of representation finite algebras is open. (For an overview
of the theory, cf \cite{GR}.)

From the preceding section, in any sequence of deformations of a
finite dimensional algebra $A$ over a field, only a finite number
can be jump deformations or effective jump deformations. Because
a representation finite algebra always has a multiplicative basis
(i.e., one where the product of two basis elements is a third or
zero) there can only be finitely many in any dimension. It
follows that any deformation of an algebra of finite
representation type which {\it effectively} changes its
structure, i.e., which is not effectively trivial, must be an
effective jump deformation. Therefore, if a representation finite
algebra $A$ is not rigid then $H^1(A,A) \ne 0$, else its group of
obstructed approximate automorphism would already be reduced to
the identity.

In view of this, call an algebra an {\it effective jump algebra}
if it admits only deformations effectively equivalent to jump
deformations. Any deformation of an effective jump algebra is
again an effective jump algebra since any finite sequence of
deformations of a finite dimensional algebra can be gathered into
one multiparameter family. It follows that the class of effective
jump algebras is open and contains the class of representation
finite algebras. Clearly an effective jump algebra can be
deformed by a finite number of effectively jump deformations to
an effectively rigid algebra. If it is not already rigid then the
preceding argument shows that it must have a non-inner
derivation.  We now ask, Is every group algebra an effective jump
algebra? By Higman's theorem a group algebra is representation
finite if and only if its  $p$-Sylow subgroups are cyclic. (For a
proof, cf. \cite{P, p.194}) Most group algebras therefore are not
representation finite but it is conceivable that all are
effective jump algebras. If so, while they can not all be
deformed to separable algebras, they all could at least be
deformed to effectively rigid ones. Note that while finite
representation type is open, it does not imply rigidity. For
example, $\Bbb F_2C_2$ has finite representation type but is not
rigid. It would be useful to know in particular if the group
algebra over $\Bbb F_2$ of the quaternion group, for which the
Donald--Flanigan conjecture fails, is an effective jump algebra.

\subhead 7. Finite abelian groups and dihedral groups
\endsubhead In this section we show that there is a global
solution to the Donald--Flanigan problem for any finite abelian
group $\Gamma$, provided that one is allowed to begin with a
suitable integral extension of $\Bbb Z$. (Donald and Flanigan
exhibited only a weak solution.) Since $\Gamma$ is a product of
groups of prime power order, it is sufficient to consider the
case of $C_r$ where $r= p^m$, a prime power. Its integral group
algebra is  $A:=\Bbb Z[x]/(x^r-1)$, which can be deformed (over
$\Bbb Z[t]$) to $A_t:=\Bbb Z[x, t]/(x^r-tx-1)$. This is not yet
separable  since we can still reduce modulo the ideal generated
by $p$ and $t$ to get an algebra with a non-trivial radical. Now
a finitely generated algebra $\Cal A$ over a commutative ring $R$
is separable if and only if $\Cal A/\frak m\Cal A$ is separable
over $R/\frak m$ for every maximal ideal $\frak m$ of $R$. Let
$d$ be the discriminant of  $x^r-tx-1$ and set $R = \Bbb Z[x, t,
1/d]$.  For $\Cal A$ now take $A_t$ with coefficients extended to
$R$, i.e., $A_t\otimes_{\Bbb Z[x,t]}R$. Equivalently, letting $S$
be the  multiplicatively closed subset consisting of the powers
of $d$, we have $\Cal A = S^{-1}A_t$. Reduction of this algebra
modulo a maximal ideal of $R$ produces an algebra over some
finite field $F$ of the form $F[x]/ (x^r-\tau x -1)$, where
$\tau$ is a non-zero element of $F$ and where $x^r-\tau x -1$ has
distinct roots because its discriminant is invertible. It is
therefore separable. Hence so is $\Cal A = S^{-1}A_t$, which is
therefore a global solution to the Donald--Flanigan  problem for
$C_r$. A problem with this solution is that the algebra is
generally not split (which here means not isomorphic to a direct
sum of copies of $R$), so there is no natural correspondence
between its blocks and those of the complex group algebra.  We
can modify the foregoing procedure to obtain a split solution.
The results may look similar but the solutions will be
inequivalent. As an example, let $\Gamma$ be the cyclic group
$C_3$ of order 3, so defining $\Cal A$ as above, the (non-split)
global deformation of its integral group algebra is given by by
setting $x^3 -tx-1 = 0$.  The discriminant of $x^3 -tx -1$ is $d
= -27 + 4t^3$ and the coefficient ring $R$ will be $\Bbb Z[t,
d^{-1}]$.   Recall that separability of an algebra $\Cal A$ over
a ring $R$ is equivalent to the existence of a separability
idempotent, i.e., an element $e = \sum x_i\otimes y_i \in \Cal
A\otimes_R\Cal A$  such that for all $a\in\Cal A$ one has  $ae =
\sum ax_i\otimes y_i = \sum x_i\otimes y_ia = ea$ and such that
$\sum x_iy_i = 1$. (There are many equivalent  definitions of the
separability of an algebra $\Cal A$ over a ring  $R$. The
existence of a separability idempotent is often the most
convenient but we shall need the homological one later.) When, as
here, $\Cal A$ is a free module over its coefficient ring $R$,
calculating the separability idempotent (which in general need
not be unique) is equivalent to solving a system of simultaneous
linear equations. If $R$ is a domain  and the system is
consistent then there is a solution in the quotient field of $R$,
but it is clear from the process that one need invert only a
single element of $R$. In the present case, $\Cal A \otimes_R
\Cal A$ is a free module over $R$ spanned by all $x^i\otimes
x^j$,  $i,j = 0, 1, 2$ and we seek a linear combination $e = \sum
c_{ij}x^i \otimes x^j$ with $\sum c_{ij}x^{i+j} = 1$ and with $xe
= \sum c_{ij}x^{i+1} \otimes x^j = \sum c_{ij}x^i \otimes x^{j+1}
= ex$, the latter condition being sufficient, since $x$ generates
$\Cal A$, to insure that $ae = ea$  for all $a \in \Cal A$. 
Replacing $x^3$ by $tx+1$ and $x^4$ by  $tx^2 + x$ gives the
equations for the $c_{ij}$.   Direct calculation then yields   
$$
\align
e = d^{-1}&[(-9+4t^3)1\otimes1 + 2t^2x\otimes x + 6t(x^2\otimes
x^2 + 1\otimes x + x\otimes 1)\\ &-4t^2(1\otimes x^2 + x^2
\otimes1) -9(x\otimes x^2 + x^2\otimes x)].
\endalign
$$ 
Note that $e$ is well-defined modulo any prime, and in
particular,   modulo $3$. If the coefficient ring had been one in
which one could divide by $3$, then at $t=0$ the separability
idempotent $e$ would reduce to 
$\frac 13(1\otimes 1 +x\otimes x^2 + x^2 \otimes x)$, i.e.,  to
the classical separability idempotent $\frac 1{\#G} \sum g\otimes
g^{-1}$ where here $\#G = 3$ and $g$ runs through the elements of
the cyclic group $C_3 = \{1, x, x^2\}$. Since the rational group
algebra  $\Bbb QC_3$  is not split the present algebra also could
not be split, else extending coefficients to include $\Bbb Q$ we
could let $t \to 0$, giving a contradiction.

To remedy this, return to the group $C_r$ with $r$ a power of a
prime $p$, let $\Phi$ be the cyclotomic polynomial for the
primitive $r$th roots of unity and take $\Cal O = \Bbb
Z[y]/\Phi(y)$. Denoting the image of $y$ by $\eta$, instead of
$x^r-tx-1$, take now the polynomial $f(x,t) = \Pi_{i=0}^{r-1}(x-
\eta^i(1+t))$. Since $f(x,0) = x^r-1$, $\Cal O[x,t]/f(x,t)$ is a
deformation of the group algebra $\Cal OC_r$. We may now proceed
exactly as before, localizing at the multiplicatively closed set
$S$ generated by the discriminant. The result is a global
solution which is split and whose blocks, therefore, are now in
natural one-to-one correspondence with those of the complex group
algebra once one fixes a morphism of $\Cal O$ into $\Bbb C$. This
will always be tacitly understood. 

To handle the dihedral groups we need one additional refinement:
a split global solution which moreover is equivariant with
respect to the automorphism sending every element of the cyclic
group to its inverse. For this we need that $f(x,t)$ be symmetric
in the sense that 
$f(x, t) =(-x)^r f(x^{-1},t)$, where $r=p^m$, as before.
Introduce a parameter $q$ which will be set equal to $1+t$ and
extend $\Cal O$ to $\Cal O[q^{-1}]$. Set  
$$ f(x,t) = \left\{
\aligned&\Pi_{i= -(r-1)/2}^{(r-1)/2}(x-\eta^i q^i)
\qquad \qquad\quad\ \text{ $p$ odd,}\\
(x- q)(x -  q^{-1}) &\Pi_{i=1}^{(r/2)-1}
(x- q^{2i}\eta^i)(x- q^{-2i} \eta^{-i}) \quad \text{$p=2.$}
\endaligned\right.
$$ 
Together we have the following.
\proclaim{Theorem 2} For every finite abelian group and every
finite dihedral group there is a finite integral extension $k$ of
$\Bbb Z$ over which there is a split global solution to the
Donald--Flanigan problem.$\qed$\endproclaim
Now consider again the example of the cyclic group of order 3.
The coefficient ring with which one starts is now $\Cal O = \Bbb 
Z[y]/(1+y+y^2): =\Bbb Z[\omega]$, where $\omega$ is the image of
$y$. If we set $s=\omega^2( q^{-1}- q^3)$ then
$f(x,t)$ now has the form $x^3-sx^2+sx-1$. Viewing $s$ as a new
parameter, this deformation is not equivalent to the preceding
since at the prime $3$ no change of variable can remove the
quadratic term from $f$. Thus we have inequivalent deformations
leading to distinct separable commutative algebras. The
discriminant now is $\Delta=(s+1)(s-3)^3$ and the separability
idempotent is
$$
\align
e = \Delta^{-1}&[(s-3)(s^3-3s^2+s+3)(1\otimes1)+2s(s-2)(s+1)(s-3)
(x\otimes x) +  \\&2s(s-3)(x^2\otimes x^2) -
s(s-2)(s+1)(s-3)(1\otimes x + x\otimes 1)
+\\&s(s-1)(s-3)(1\otimes x^2 + x^2\otimes 1) - 
(2s-3)(s+1)(s-3)(x\otimes x^2 + x^2\otimes x)].
\endalign
$$ 
This is again well-defined modulo $3$ and reduces to
$\frac 13(1\otimes 1 + x\otimes x^2 + x^2\otimes x)$ at $s=0$.
Recalling that $ q = 1+t$, this is equivalent to setting
$t=0$. (While $\Cal O$ is not separable over $\Bbb Z$ since
modulo $3$ it contains a central nilpotent element the deformed
algebra is separable over $\Cal O$, which is all we want.) This
new solution is split when $s$ is replaced by
$\omega^2( q^{-1}- q^3)$. Although it is inequivalent to the
preceding one, this does {\it not} imply that over $\Bbb F_3$ the
algebra defined by $x^3 = 1$ lies on the intersection of two
components of the variety of algebras of dimension 3. For over an
algebraically closed field there is only one separable algebra of
dimension 3, namely, the sum of three copies of the coefficient
field, so the deformations lie on the same component. But only
for the second (split) one is there a fixed correspondence
between its summands and those of $\Bbb CC_3$. 

We will show, in particular, that i) if the Donald--Flanigan
problem has a solution (in the weak sense, i.e., at a fixed
prime) for some group $K$ then it also has a solution for the
wreath products $K\wr \Gamma, K\wr S_n$ of $K$ with any abelian
or symmetric group and ii) if there is a split global solution
for $K$, i.e., one which is a sum of matrix algebras, then the
same is true for these wreath products. This will imply, in
particular, that the problem has a global solution for the the
Weyl groups of type $B_n$ and by  minor modifications for the
dihedral groups and groups of type $D_n$. For the groups $B_n$ 
and $D_n$ it will not be necessary to
adjoin any roots of unity to $\Bbb Z$, so we recover the fact
that all their irreducible complex representations also are
rational. If instead of a split global solution for $K$ we have
only a central separable one, then as mentioned in the
introduction we shall need some conditions on the center, $\Cal
Z$ of the deformed algebra. Suppose that the new coefficient ring
is $R$. To insure that $\Cal Z$ can be split by a separable,
faithfully projective extension of $R$ we shall have to assume,
that $\Cal Z$ is faithfully projective and of constant rank over
$R$. 

\subhead 8. Some smash and wreath products \endsubhead
We will need to know the structure of a smash product $A\#kG$
where $G$ is a finite group and $A$ a $k$-algebra of the form
$k^n = k\oplus \cdots \oplus k\, \text{($n$ times)}$. The only
way that a group $G$ can act as $k$-algebra automorphisms of such
an algebra $A$ is by permutation of the summands, in effect, by
permutation of the indices $1,\dots,n$. Denoting the unit element
of the $i$th summand of $A$ by $e(i)$, the smash product is
spanned by the elements $e(i)\sigma, \sigma \in G$. These
multiply by the rule $\sigma e(i) = e(\sigma i)\sigma$, so
$e(i)\sigma e(j)\tau = \delta_{i,\sigma j}e(i)\sigma\tau$. It
follows that $A\#kG$ is a direct sum of subalgebras corresponding
to the orbits of $G$. Fix some $i\in \{1,\dots,,n\}$, denote the
orbit of $i$ by $Gi$, the order of this orbit by $m_i$, and the
isotropy group of $i$ by $G_i$. If
$\{\bar\sigma,\bar\tau,\dots\}$ are representatives of the cosets
$G/G_i$ then the various $e(\bar\sigma i)\bar\sigma^{-1}
(\bar\tau)$ span a subalgebra of $A\#kG$ isomorphic to
$M_{m_i}(k)$. For writing  $e(\bar\sigma i)(\bar\tau)^{-1}=
E_{\bar\sigma,\bar\tau}$, these multiply just like matrix units,
i.e., $E_{\bar\mu,\bar\nu}E_{\bar\sigma,\bar\tau} =
\delta_{\bar\nu,\bar\sigma}E_{\bar\mu,\bar\tau}$. 
\proclaim{Theorem 3} If $A = k^n$ and $G$ is a finite group
operating by permutation of the summands then
$$
A\#kG \cong \bigoplus_i M_{m_i}(k)\otimes kG_i
$$ where the sum is over representatives $i$ of the distinct
orbits of $G$ in $\{1,\dots,n\}$, $m_i$ is the order of the orbit
$Gi$ of $i$, and $G_i$ is the isotropy group of $i$.
\endproclaim
\demo{Proof} Since $A\#G$ is the direct sum of subalgebras
corresponding to the orbits of $G$, what we must show is that the
subalgebra corresponding to the orbit of $i$ is isomorphic to
$M_{m_i}\otimes kG_i$. Letting $\bar\sigma,
\bar\tau, \dots$ be as before representatives of the cosets
$G/G_i$, the $E_{\bar\sigma,\bar\tau} = e(\bar\sigma
i)\bar\sigma^{-1}\bar\tau$ span an algebra isomorphic to
$M_{m_i}(k)$. This contains, in particular, the permutation
matrices, so if we have any permutation $\pi$ of $\{1,\dots,n\}$
there is a unique matrix $T_{\pi} \in M_{m_i}(k)$ giving the same
permutation. Thus the operator $T_{\pi}^{-1}\pi$ acts as the
identity on $\{1,\dots,n\}$ and hence commutes with the operation
of all elements of $M_{m_i}(k)$. Therefore $T_{\pi}^{-1}\pi
T_{\rho}^{-1}\rho = T_{\rho}^{-1}T_{\pi}^{-1}\pi\rho =
T_{\pi\rho}^{-1}\pi\rho$. The map $\pi \to T_{\pi}^{-1}\pi$ is
therefore a group morphism. Since the elements $E_{\bar\sigma,
\bar\tau}\rho$ with $\rho \in G_i$ form a $k$-basis for $A$, this
proves the assertion. $\qed$
\enddemo
Note here that the elements of the isotropy group do not
themselves commute with the elements of $M_{m_i}$. The second
tensor factor in the statement of theorem is not the group
algebra of the isotropy group $G_i$ itself but the group algebra
of the isomorphic group of all $T_{\pi}^{-1}\pi$. An analogous
argument is used in the next two theorems. One consequence of the
above theorem is that if we have an algebra of the form $A =
k^n\#kG$ such that the Donald--Flanigan problem is solvable for
all the isotropy groups $G_i$ then it is solvable for all of $A$.

In addition to the foregoing we need to know the structure of
$A\wr G$ when $A$ is a direct sum of central separable algebras
over some coefficient ring $k$. In this section we compute this
when $A$ is itself central separable and find that $A\wr G$ is
canonically isomorphic to $A^{\otimes n}\otimes kG$. This makes
it possible to extend any deformation of $kG$ to one of $A\wr G$.
Again, while $kG$ is contained in a natural way in both $A\wr G$
and $A\otimes G$, the isomorphism, while canonical, does not
carry one embedding to the other. The result is extended to the
general case in the next section.

Denote the group of invertible elements of $A$ by $A^*$, its
group of automorphisms by $\operatorname{Aut}A$, and the center
of $A$ by $\Cal Z$. There is a group morphism $A^* \to
\operatorname{Aut}A$ sending $a\in A^*$ to the inner automorphism
$\operatorname{in}_a$ defined by $\operatorname{in}_ax = 
axa^{-1}$; the kernel of this morphism is $\Cal Z^*$. For a
central simple algebra over a field the Skolem--Noether theorem
asserts that the morphism is onto, but this need not be so for an
arbitrary central separable algebra. Action of a group $G$ on $A$
is the same as a group morphism $f:G\to \operatorname{Aut}A$. The
morphism $f$ ``factors through'' the canonical morphism $A^* \to
\operatorname{Aut}A$ if there is a morphism $G \to A^*$ whose
composite with the canonical morphism $A^* \to
\operatorname{Aut}A$ is $f$. This is stronger than asserting that
the image of every individual $\sigma \in G$ is inner (which
would always be the case for a central simple algebra). For if
$\sigma \in G$ then the choice of an element $a_{\sigma} \in A$
conjugation by which has the same effect as $\sigma$ is
determined only up to an element of $\Cal Z^*$. It may not be
possible to make the choices in a consistent way to produce a
morphism $G \to A^*$, as the following trivial example shows. Let
$A = M_2(k)$ and set $a = \pmatrix 1 & 0\\0 & -1 \endpmatrix,
\quad b = \pmatrix 0 & 1\\1 & 0 \endpmatrix$. Then
$\operatorname{in}_a$ and  $\operatorname{in}_b$ commute and both
have square equal to the identity, so they define an operation of
$C_2 \times C_2$ on $M_2(k)$. But $ab = -ba$ and multiplying $a,
b$, or both by an invertible central element can not change this.
Therefore (if the characteristic is not 2) the morphism of
$C_2\times C_2$ into $\operatorname{Aut} A$ does not factor
through $A^*$. The problem is homological. For choosing an
$a_{\sigma}$ for every $\sigma$ defines an element of $H^2(G,\Cal
Z^*)$. This class does not depend on the choices and may not be
trivial. (If $k$ contains a square root of -1, which we may
denote by $i$,  then $\operatorname{in}_{ia} =
\operatorname{in}_a, \operatorname{in}_{ib}  =
\operatorname{in}_b$, and $ia, ib$ generate a group 
isomorphic to the quaternion group.) 
\proclaim{Lemma A} Supppose that a group $G$ operates as a group
of automorphisms on a $k$-algebra $A$ and that the associated
morphism $G \to \operatorname{Aut}A$ factors through $A^*$. Then
$A\#kG \cong A\otimes kG$. 
\endproclaim
\demo{Proof} To avoid confusion, denote the image of $x \in A$
under the operation of $\sigma \in G$ by $x^{\sigma}$. By
hypothesis, for each $\sigma \in G$ there exists an invertible
element $a_{\sigma}$ such that $a_{\sigma}xa_{\sigma}^{-1} =
x^{\sigma}$ and $a_{\sigma}a_{\tau} = a_{\sigma\tau}$ all
$\sigma,\tau \in G$. It follows that the elements 
$a_{\sigma}^{-1}\sigma \in A\#kG$ commute with all elements of
$A$ and form a multiplicative subgroup of $G'$ of $A\#kG$
isomorphic to $G$. But then $A\#kG$ = $A\otimes kG' \cong
A\otimes kG$. For if $\{a_i\}$ is a $k$-basis for $A$ then every
element of $A\#kG$ can obviously be written uniquely as a linear
combination with coefficients in $k$ of the elements $a_i\otimes
a_{\sigma}^{-1}\sigma$. $\qed$ \enddemo

Suppose now that $A$ is a central separable $k$-algebra. Every
$a\in A$ then has a reduced trace $\operatorname{Trd}(a) \in k$;
when $A = M_n(k)$ this is just the usual trace of a matrix (cf.
\cite{KO, pp. 90, 110}). (The ``ordinary'' trace of $a$ comes
from the characteristic polynomial of left multiplication by $a$
in the algebra. This can vanish even for a primitive idempotent.
Consider e.g. $2\times 2$ matrices in characteristic $2$. Here
the trace of left multiplication by $E_{11}$ is zero but the
reduced trace is 1.)  Now a central separable $k$-algebra $A$
always contains $k$ as a $k$-module direct summand so the map
$a\mapsto \operatorname{Trd}a$ may be viewed as $k$-module
endomorphism of $A$. Since for central separable $A$ there is an
isomorphism $A\otimes A^{\operatorname{op}} \to
\operatorname{End}_kA$ sending $x\otimes y$ to the endomorphism
carrying $a\in A$ to $xay$, it follows that there exists a unique
element $T = \sum x_i \otimes y_i \in A\otimes A$ such that
$\operatorname{Trd}a = \sum x_iay_i$. Viewing this ``switch''
element $T$ as an element of $A\otimes A$, it has the properties
i) $T^2 = 1\otimes1$, the unit element of $A\otimes A$ and ii)
$T(a\otimes b)T = b\otimes a$ for all $a, b\in A$. (This result,
communicated by A. Fr\"olich to Knus and Ojanguren is credited by
them \cite{KO, p.112} to O. Goldman who with M. Auslander
developed the
present concept of separable algebra, \cite{AG}.)  
For $A = M_n(k)$ one has $T = \sum_{i,j} e_{ij} \otimes
e_{ji}$. It follows from the preceding Lemma that if $C_2$
operates on $A\otimes A$ by interchange of the tensor factors,
then $A\wr C_2 = A^{\otimes 2}\#kC_2$ is canonically isomorphic
to $A^{\otimes 2}\otimes kC_2$. The isomorphism is canonical
because there is a canonical choice of an element in $A\otimes A$
conjugation by which interchanges the tensor factors, namely $T$.
Note, incidentally, that $T$ is symmetric since $T = T^3 =
T(T)T$, i.e., $T$ with its tensor factors interchanged.

\proclaim{Theorem 4} Let $A$ be a central separable $k$-algebra
and let $S_n$ operate on $A^{\otimes n}$ by permutation of the
tensor factors. Then there is a canonical isomorphism $A\wr S_n
\cong A^{\otimes n} \otimes kS_n$. 
\endproclaim
\demo{Proof} The group $S_n$ is generated by elements $T_{12},
T_{23}, \dots, T_{n-1,n}$ corresponding to the transpositions
$(12), (23), \dots, (n-1,n)$  subject only to the conditions that
$T_{i,i+1}$ and $T_{j,j+1}$ commute for $|i-j| > 1$, that
$T_{i,i+1}^2 = 1$, and the braid relation $T_{i, i+1}T_{i+1,i+2}
T_{i,i+1} = T_{i+1, i+2}T_{i,i+1}T_{i+1,i+2}$; this is the
``Artin presentation'' of $S_n$. Now let $T_{i,i+1}$ denote
the switch element $T$ operating in places $i$ and $i+1$ of
$A^{\otimes n}$. It will be sufficient to show that these satisfy
the relations of the Artin presentation. Since $T^2 = 1$ the
first relations are clearly satisfied. For the braid relations it
is only necessary to
consider the case $n=3$ where we must prove that
$T_{12}T_{23}T_{12} =
T_{23}T_{12}T_{23}$. If $T = \sum x_i \otimes y_i$ then the left
side is $\sum (x_i\otimes y_i \otimes 1)(1\otimes x_j\otimes
y_j)(x_k \otimes y_k \otimes 1) = \sum x_ix_k \otimes y_ix_jy_k
\otimes y_j$. However, $\sum x_ix_k \otimes y_ix_jy_k =
T(1\otimes x_j)T = x_j\otimes 1$, so this is just $\sum x_j
\otimes 1 \otimes y_j$ (which would naturally be denoted
$T_{13}$). One checks similarly that the right side gives the
same thing. $\qed$
\enddemo
It was convenient in the preceding to deal with $S_n$ but it is
clear that nothing would change if instead of $S_n$ itself we had
a group $G$ which operated by permutation of the tensor factors,
i.e., through a morphism $G \to S_n$.  This gives the following.
\proclaim{Corollary} If $A$ is central separable over $k$ and a
group $G$ operates on $A^{\otimes n}$ through a morphism $G \to
S_n$ then $A\wr G$ is canonically isomorphic to $A^{\otimes n}
\otimes kG$. $\qed$
\endproclaim

It follows that any deformation of $kH$ induces one of $A\wr H$.

\subhead 9. Wreath product with a sum of central separable
algebras \endsubhead
Suppose now that the $k$-algebra $A$ is a direct sum of central
separable algebras, $A = A_1 \oplus \dots \oplus A_r$ with
respective unit elements $e_1,\dots,e_r$ and let $\Cal Z$ denote
the commutative algebra which these generate.  Extending the
result of the previous section, we wish to determine the
structure of $A\wr S_n$ and more generally of $A\wr G$ where $G$
is a group operating on $A^{\otimes n}$ by permutation of the
tensor factors. For any multiindex $I = (i_1,\dots,i_n)$ set
$A(I) = A_{i_1}\otimes A_{i_2} \otimes \cdots \otimes A_{i_n}$
and denote its unit element $e_{i_1}\otimes e_{i_2}\otimes \cdots
\otimes e_{i_n}$ by $e(I)$. Then $A^{\otimes n}$ is the direct
sum of the $A(I)$, which are all central separable, one has
$e(I)e(J) = \delta_{IJ}e(I)$, and the distinct $A(I)$ are
mutually orthogonal. The wreath product $A\wr G$ is spanned by
elements of the form $a(I)\otimes \sigma, a(I) \in A(I), \sigma
\in G$, which we write simply as $a(I)\sigma$. One has
$e(I)\sigma e(J)\tau = e(I)e(\sigma J)\sigma\tau$, where 
$\sigma (j_1,\dots,j_n) = (j_{\sigma^{-1}1},\dots,
j_{\sigma^{-1}n})$. We write $(a_1,\dots,a_n)$ for $a_1\otimes
\cdots \otimes a_n \in A^{\otimes n}$ where each $a_i$ is in $A$,
and if $\alpha = (a_1,\dots,a_n)$ then $\sigma\alpha =
\alpha^{\sigma} \sigma$, where $(a_1,\dots,a_n)^{\sigma} =
(a_{\sigma^{-1}1}, \dots, a_{\sigma^{-1}n})$.   

Set $B = A\wr G$, choose a multiindex $I$, and let $B(I)$ be the
subalgebra  generated by $A(I)$ and the elements of $G$. It is
spanned by all $\alpha^{\sigma}\tau$.  One has $B(I) =
B(\sigma I)$ for all $\sigma \in G$, so these subalgebras are
really indexed by the orbits of $G$ in the set of multiindices,
and $B$ is their direct sum. Determining the structure of $B$ is
thus the same as determining that of the various $B(I)$. If there
is a global solution to the Donald--Flanigan problem for each of
the ``orbital'' subalgebras $B(I)$ then the same will be true for
$A\wr kG$, and similarly if there is only a weak solution for
each of the subalgebras then the same will be hold for $A\wr kG$.
Fixing the multiindex $I$, let $G_I$ be the isotropy group of
$I$. Now choose fixed representatives $\bar\sigma, \bar\tau,
\dots$ for the cosets $\sigma G_I,\tau G_I,\dots$. The elements 
$e(\bar\sigma I) \bar\sigma\bar\tau^{-1}$ then multiply precisely
like the matrix units $E_{\bar\sigma,\bar\tau}$ (whose rows and
columns are indexed by the chosen coset representatives), namely
$(e(\bar\rho I)\bar\rho\bar\mu^{-1})(e(\bar\sigma I)
\bar\sigma\bar\tau^{-1}) = \delta_{\bar\mu\bar\sigma}e(\bar\rho
I)\bar\rho\bar\tau^{-1}$. 

The choice of coset representatives $\{\bar\mu,
\bar\sigma,\dots\}$ fixes particular isomorphisms $A(I) \to
A(\bar\sigma I)$ and hence also isomorphisms $A(\bar\sigma I) \to
A(\bar\mu I)$ for all $\bar\sigma$ and $\bar\mu$. We can write
the elements of $A(\bar\sigma I)$ in the form $a^{\bar\sigma}$
with $a \in A(I)$. Since
$(a^{\bar\mu}\bar\mu\bar\sigma^{-1})(b^{\bar\sigma}\bar\sigma^{-1
}\bar\tau) = (ab)^{\bar\mu}\bar\mu\bar\tau^{-1}$ the elements of
the form $a^{\bar\mu}\bar\mu\bar\sigma^{-1}$ span a subalgebra
$\Cal M(I)$ of $A\wr G$. Let the index of $G_I$ in $G$ be $m$
(the cardinality of the set of representatives). Then the
subalgebra $\Cal M(I)$ is isomorphic to $M_m(A(I))$ under the map
$a^{\bar\sigma}\bar\sigma\bar\tau^{-1} \mapsto
aE_{\bar\sigma,\bar\tau}$. We will show that $B(I) \cong
M_m(A(I)) \otimes kG_I$, but as in the previous section, the
tensor factor $kG_I$ here is not the usual subalgebra $kG_I$ of
$A\wr G$. With the foregoing notations and writing $\alpha^{\pi}$
for $\pi(\alpha), \alpha \in A(I), \pi \in G$ we have the
following.
\proclaim{Lemma B} For every $\pi \in G_I$ there is a canonical
choice of $T_{\pi} \in A(I)$ such that $T_{\pi}\alpha
T_{\pi}^{-1} =
\alpha^{\pi}$. \endproclaim
\demo{Proof} Suppose that $A(I) = A(i_1,\dots, i_n) = A_{i_1}
\otimes \cdots \otimes A_{i_n}$. For convenience, reorder the
tensor factors, putting the indices in increasing order, so that
the tensor product has the form $A_1^{\otimes n_1} \otimes
\cdots \otimes A_r^{\otimes n_r}$ where $n_1 + \dots +n_r = n$.
Since $\pi I = I$, after the reordering $\pi$ permutes the first
$n_1$ factors amongst themselves, the next $n_2$ factors amongst
themselves, etc.. From the preceding section we know that there
is a canonical $T^{(1)}_{\pi}\in A_1^{\otimes n_1}$ such that for
any $\alpha_1 \in A_1^{\otimes n_1}$ we have
$T^{(1)}_{\pi}\alpha_1 (T^{(1)}_{\pi})^{-1} = (\alpha_1)^{\pi}$,
and similarly for the other $A_i^{\otimes n_i}$. Then
$T^{(1)}_{\pi} \otimes \cdots \otimes T^{(r)}_{\pi}$ with factors
returned to their original order is the desired $T_{\pi}$. $\qed$
\enddemo

Now set $\upsilon(\pi) = \sum \bar\mu T^{-1}_{\pi}\pi \bar
\mu^{-1}$, where the sum runs over all coset representatives
$\bar \mu$ of $G_I$. Since $\bar\mu T^{-1}_{\pi}\pi \bar \mu^{-1}
= (T_{\pi}^{-1})^{\bar\mu}\bar\mu\pi\bar\mu^{-1} 
\in A(\bar \mu I)\bar\mu\pi\bar \mu^{-1}$ the terms for different
$\bar \mu$ are mututally orthogonal, from which it is easy to see
that $\upsilon(\pi)\upsilon(\pi') = \upsilon(\pi\pi')$ for $\pi,
\pi' \in G_I$. Since $\upsilon(1) = \sum \bar\mu e(I)
\bar\mu^{-1} =\sum e(\bar \mu I)$ is just the identity element of
$B(I) = \oplus A(\bar\mu I)$, the $\upsilon(\pi)$ form a
multiplicative subgroup of $B(I)$ isomorphic to $G_I$. Moreover,
if $\alpha \in A(I)$ then $\upsilon(\pi)$ also commutes with all
$a^{\bar \sigma}\bar\sigma \bar \tau^{-1}$.

\proclaim{Lemma C} There is a canonical isomorphism $B(I) \cong
M(I) \otimes kG_I$ given by
$a^{\bar\sigma}\bar\sigma\bar\tau^{-1}\upsilon(\pi) \mapsto
aE_{\bar\sigma,\bar\tau}\otimes\pi$. 
\endproclaim
\demo{Proof} Observe that $a^{\bar\sigma}\bar\sigma\bar\tau$ is
orthogonal to all the summands of $\upsilon(\pi) = \sum\bar\mu
T_{\pi}^{-1}\pi\bar\mu^{-1}$ except that with $\bar\mu =
\bar\tau$, so $a^{\bar\sigma}\bar\sigma\bar\tau^{-1}\upsilon(\pi)
= (aT_{\pi}^{-1})^{\bar\sigma}\bar\sigma\pi\bar\tau^{-1} =
(aT_{\pi}^{-1})^{\bar\sigma}\bar\sigma(\bar\tau\pi^{-1})^{-1}$.
As $a$ runs through $A(I)$, $(aT_{\pi}^{-1})^{\bar\sigma}$ runs
through $A(\bar\sigma I)$. As $\pi$ runs through $G_I$ so does
$\pi^{-1}$ and $\bar\tau\pi^{-1}$ runs through the coset of
$\bar\tau$. Consequently, the $\bar\tau\pi^{-1}$ exhaust $G$ and
therefore so do their inverses, and hence, for any fixed
$\bar\sigma$ so do the $\bar\sigma\pi\bar\tau^{-1}$. So for fixed
$\bar\sigma$, the set of all
$a^{\bar\sigma}\bar\sigma\bar\tau^{-1}\upsilon(\pi)$ with varying
$\bar\tau, \pi$ is precisely $B(I)$, from which the result
follows. $\qed$
\enddemo
We therefore have the following.
\proclaim{Theorem 5} Let $A = A_1\oplus\cdots\oplus A_r$ be a sum
of central separable $k$-algebras $A_i$, and suppose that a group
$G$ operates on $A^{\otimes n}$ by permutation of the tensor
factors. For every multiindex $I= (i_1,\dots,i_n)$ set $A(I) =
A_{i_1}\otimes\cdots\otimes A_{i_n}$, let $G_I$ be the isotropy
group of $I$, and let $\{\bar I, \bar J, \dots\}$ be a set of
representatives of the orbits of $G$ in the set of multiindices.
Set $m_I = (G:G_I)$, the index of $G_I$ in $G$. Then $A\wr G$ is
canonically isomorphic to $\bigoplus M_{m_{\bar I}}(A_{\bar I})
\otimes kG_{\bar I}$ where the sum runs over the orbits of $G$ in
the set of multiindices. $\qed$
\endproclaim

There is an overlap between this theorem and Theorem 3, namely,
the case where each $A_i$ is just $k$ itself and one has a wreath
product, but Theorem 3 allows more general smash products. As an
immediate corollary of Theorem 5 one gets explicitly the
structure of the complex group algebra (and hence the irreducible
complex representations) of a wreath product of groups $H\wr G$.
More generally, if $k$ is a splitting ring for both $H$ and all
the $G_I$, i.e., one such that both $kH$ and the $kG_I$ are sums
of central separable algebras, then the theorem implies that $k$
is also a splitting ring for $H\wr G$ and it gives the structure
of $k(H\wr G)$. Consider, for example, $B_n = C_2\wr S_n$. Here
$\Bbb Q$ is a splitting field for $C_2$. Writing $C_2 = \{1, a\}$
and setting $e = (1+a)/2, f = (1-a)/2$, we have $\Bbb QC_2 = \Bbb
Qe \oplus \Bbb Qf$, a direct sum of two copies of $\Bbb Q$. The
idempotents $e(I)$ are here just tensor products of length $n$ of
$e$'s and $f$'s. If $e$ occurs $m$ times and $f$ occurs $n-m$
times then the order of its orbit is $\binom nm$ and the
corresponding isotropy group is isomorphic to $S_m\times 
S_{n-m}$. Since $\Bbb Q$ is also a splitting field for the
symmetric group, it is a splitting field for $B_n$, so the
irreducible complex representations of $B_n$ are all rational.
The theorem then gives
$$
\Bbb QB_n = \bigoplus_{m=0}^n M_{\binom nm}(\Bbb Q) \otimes \Bbb
QS_m \otimes \Bbb QS_{n-m}.
$$
This includes the classic result that the representations of
$B_n$ are indexed by all {\it ordered} pairs consisting of a
representation of $S_m$ and a representation of $S_{n-m}$. For if
$\lambda, \mu$ are partitions of $m$ and $n-m$, respectively, and
if $S_m(\lambda)$ is the block (=matrix algebra summand) of $\Bbb
QS_m$ corresponding to $\lambda$ and $S_{n-m}(\mu)$ the block of
$\Bbb QS_{n-m}$ corresponding to $\mu$, then the blocks of $\Bbb
QB_n$ are of the form $M_{\binom nm}\otimes
S_m(\lambda)\otimes   S_{n-m}(\mu)$. The representations $B_n$
are thus indexed by the pairs 
$(\lambda, \mu)$. Reversing their order gives a different
representation. By contrast (cf. \S11), the representations of
$D_n$ may be viewed as indexed by {\it unordered} pairs of
representations of $S_m$ and $S_{n-m}$ (or by pairs in which
$m\le n-m$).

\subhead 10. First applications\endsubhead 
The most immediate application of the results of the preceding
section is to the case where we have a group $H$ for which there
is a split global solution to the  Donald--Flanigan problem and
we seek one for a group of the form $H\wr G$. So suppose that $k$
is some subring of algebraic integers in $\Bbb C$ and that $kH$
has been deformed into a direct sum $A = M_{i_1}(R) \oplus\cdots
\oplus M_{i_r}(R)$ of matrix algebras. Here the new coefficient
ring $R$ is the localization of a power series ring
$k[[t,u,\dots,v]]$ at a single element $f$. Since $R$ is a
domain, the inclusion of $k$ into $\Bbb C$ can be extended to a
monomorphism of $R$ into $\Bbb C$. If we have been careful with
$R$ to require, e.g., that it is projectively faithful over $k$,
then we can still reduce modulo every prime $p$ of $k$ (see
below), so there is a natural correspondence between the blocks
of $\Bbb CH$ and those we get by reducing $A$ modulo $p$.  By
assumption there is given some morphism $G \to S_n$, and through
this $G$ operates on $H^{\otimes n}$ by permutation of the tensor
factors.  Keeping the notations, the results of the preceding
section immediately yield the following. 

\proclaim{Theorem 6} If there is a split global solution for $RH$
and if in $H\wr G$ for each isotropy group $G_I$ there is a split
global solution to the Donald--Flanigan problem for $RG_I$ then
there is one for $H\wr G$ over $k$. (If there is only a weak
solution for $RG_I$ then there is a weak solution for $H\wr G$
and similarly for local solutions.)$\qed$
\endproclaim

There are several important cases in which the conditions on $G$
will be fulfilled. First, if $G$ is abelian then so is every
isotropy group; we have already constructed split global
solutions for finite abelian groups. Second, suppose that $G$ is
all of $S_n$. The isotropy group of an index of the form
$(i_1^{n_1}, \dots, i_r^{n_r})$, where $i^m$ stands for $i$
repeated $m$ times and $n_1\dots +n_r = n$ is $S_{n_1} \times
\cdots \times S_{n_r}$. Since we have split global solutions for
the symmetric groups we also have a split global solutions for
these. 

\proclaim{Theorem 7} There is a split global solution to the 
Donald--Flanigan problem for any wreath product $H\wr G$ where
$H$ is finite abelian and $G$ is either finite abelian or
$G=S_n$. $\qed$
\endproclaim

The Weyl groups of type $B_n$ are a special case. Since both for
$C_2$ and $S_n$ we can start with the coefficient ring $\Bbb Z$
we do not have to adjoin any roots of unity. 
\proclaim{Theorem 8} There is a split global solution to the
Donald--Flanigan problem for the Weyl groups of type $B_n$
starting with $\Bbb Z$. For the coefficient ring after
deformation one can take $R = \Bbb Z_{q,n}$ since there is a
split global deformation of $S_n$ defined over this $R$. It
follows that all complex irreducible representations of the Weyl
groups of type $B_n$ are rational.
$\qed$
\endproclaim

To handle dihedral groups in this context, consider $S_2$ acting
on an abelian group $K$ by sending every element of the latter to
its inverse. We can form the semidirect product $H\rtimes S_2$;
for $H$ the cyclic group $C_n$ this is just the dihedral group
$\Cal D_n$. Let $A$ be a split equivariant solution to the
Donald--Flanigan problem for $H$. (Here, if $e$ is the exponent
of $H$ then we must work over $\bar \Bbb Z = \Bbb Z[y]/\Phi(y)$,
where $\Phi$ is the cyclotomic polynomial for the primitive  $e$-
th roots of unity.) This $A$ is then just a direct sum of copies
of whatever coefficient ring $k$ has been introduced, and $S_2$
acts as automorphisms of the sum. However, since the automorphism
group of each summand is reduced to the identity, $S_2$ must
either leave an individual summand fixed or exchange them in
pairs. The equivariant deformation has thus produced a direct sum
of copies of $kS_2$ plus copies of $k \oplus k$ on which $S_2$
operates by interchange of the summands. As a special case of
what we have seen before, the latter yields $2 \times 2$ matrix
algebras over $k$, while $kS_2$ itself has a global deformation
to a sum of two copies of $k$ (extended by the deformation
parameter).  Therefore we have

\proclaim{Theorem 9} There is a split global solution to the
Donald--Flanigan problem for any semidirect product $H\rtimes
S_2$ where $H$ is a finite abelian group. $\qed$
\endproclaim

Notice that the argument implies the familiar fact that all
complex irreducible representations of a dihedral group $\Cal
D_m$ have dimension either one or two and and are defined over
$\Bbb Q(\eta)$ where $\eta$ is a primitive $m$-th root of unity 
(but in particular cases may be defined over a smaller field).
Since every subgroup of a dihedral group is again either dihedral
or abelian we have the following.

\proclaim{Theorem 10} There is a split global solution to the
Donald--Flanigan problem for any group of the form $H\wr\Cal D_m$
with $H$ abelian. $\qed$
\endproclaim
Suppose now that $H$ is a group for which the Donald--Flanigan
problem is solvable over the coefficient ring $k$ but where the
solution is a non-split separable algebra $A$ over the new
coefficient ring $R$. Now $A$ is Azumaya (central separable) over
its center $\Cal Z$, which in turn is separable over $R$ (cf.,
e.g., \cite{DI}, p. 55). We will be able to use the preceding
results if we can ``split'' $\Cal Z$ by a suitable extension $T$
of the coefficient ring $R$, i.e., reduce it to a direct sum of
copies of the new coefficient ring.  For this we can apply the
following basic result: For a commutative $R$-algebra $S$ the
following properties are equivalent: 1) $S$ is separable,
faithfully projective (i.e., faithful and projective) and of
constant rank over $R$, 2) there is a faithfully projective and
separable $R$-algebra $T$ splitting $S$, i.e., such that
$S\otimes_R T \cong T\oplus T \oplus \cdots \oplus T \text{ ($n$
times)} = T^n$ and 3) there is a faithfully flat and separable
$R$-algebra $T$ such that $S\otimes_R T \cong T^n$ (cf.
\cite{KO}, Th\'eor\`eme 4.7, p. 88). Projectivity and
faithfulness of T will insure, in particular, that i) if $p$ is a
rational prime then the sequence $0\to pA\otimes T \to A\otimes T
\to (A/pA)\otimes T \to 0$ is exact and ii) the map of $A$ into
$A\otimes T$ given by $a \mapsto a\otimes 1$ is an inclusion, and
similarly for $A/pA$ in place of $A$.  Also, if we have an
inclusion $R\hookrightarrow \Bbb C$ then we still have an
inclusion $T = R\otimes T \hookrightarrow \Bbb C \otimes T$
(where $\Bbb C$ is regarded as an $R$-module by the first
inclusion. 

If $k$ is a direct sum or product of subrings $k_i$ then the
Donald--Flanigan problem for a group $G$ will be solvable if it
is solvable for each $k_iG$. To avoid the complications of
infinite sums or products, we may restrict attention to the case
where $k$ has only finitely many idempotents. This holds, in
particular, if $k$ is noetherian, which is always the case here.
With this, we can reduce to the case where the only non-zero
idempotent in $k$ is its identity. A projective module over $k$
must then have constant rank. If $k$ has only finitely many
idempotents, so $k = k_1\oplus \cdots \oplus k_r$, where each
summand has no idempotent but the identity, then any $k$-algebra
$A$ is similarly a sum of $k_i$-algebras $A_i$. If each of these
can be deformed to a separable $k_i$-algebra then $A$ can be
deformed to a separable $k$-algebra. 

Deformation can not enlarge the center of an algebra but it can
diminish it. We will say that a deformation is {\it
centrality preserving} if elements which were central before
deformation continue to be central. (It need not be center
preserving because in general the center is deformed.) This must
be the case, in particular, for a solution to the 
Donald--Flanigan problem for $kG$ which resembles the complex
group algebra.  For any finite group $G$ and
coefficient ring $k$, the center of $kG$ is a free module over
$k$; a typical basis element is just the sum of the elements of
some conjugacy class. The center of a centrality preserving
deformation of a group algebra thus continues to be a free module
over the new coefficient ring $R$, so condition 1) above will
hold without the assumption that $k$ is noetherian. (Of course it
holds if the original group is commutative and the deformed group
algebra remains so, but there are important examples of
deformations where the commutativity is lost.) 

\subhead 11. Weyl groups of type $D_n$ \endsubhead 
To describe the Weyl groups of type $D_n$, write the elements of
$B_n$ in the form $(c_1,\dots,c_n)\sigma$ with $c_i \in C_2,
\sigma \in S_n$; then $D_n$ is the subgroup consisting of all
elements in which an even number of the $c_i$ are equal to the
unit element of $C_2$. Equivalently, the product of the $c_i$ is
equal to the unit element. So viewed, $D_n$ is an example of a
class of semidirect products which are subgroups of wreath
products arising as follows. Let $H$ be any abelian group (not
necessarily finite) and consider the morphism $\mu: H^n = H
\times \cdots \times H\text{ ($n$ times )} \to H$ sending 
$(c_1,\dots,c_n)$ to the product $c_1c_2\dots c_n$. If $K$ is any
subgroup of $H$ then $\mu^{-1}K$ is stable under the operation of
$S_n$ on the factors of $H^n$ so one can form the semidirect
product $(\mu^{-1}K)\rtimes S_n$. For $n=2$, $H$ cyclic and $K$
reduced to the unit element one gets the dihedral groups. For $H
= C_2$ and $K$ the unit subgroup one gets the groups $D_n$.
Writing simply $C_2^n$ for $C_2^{\times n}$, here $\mu^{-1}K
\cong C_2^{n-1}$ so the Weyl group $D_n$ is a semidirect product
$C_2^{n-1}\rtimes S_n$. (One has  $D_2 =S_2, D_3 =S_4$.)

It will be convenient to deal with $D_{n+1} = C_2^n\rtimes
S_{n+1}$ rather than $D_n$. The deformation of $C_2^n$ we
employed in constructing a solution for $B_n$ was equivariant
under the operation of $S_n$, but in $D_{n+1}$ this group is
acted upon by $S_{n+1}$ and the deformation ceases to be
equivariant. The original operations of $S_{n+1}$ consequently
cease to be automorphisms, but the operations can also be so
deformed that $S_{n+1}$ continues to operate as automorphisms of
the deformed algebra. This would normally be a difficult
homological problem but fortunately here there is a natural
solution (although it exhibits a rather unnatural phenomenon at
the prime $p = 2$).  To exhibit it, however, we first examine the
rational group algebra of $D_{n+1}$. 

To describe the operation of $S_{n+1}$ on $C^n$, view $S_{n+1}$
as permutations of the set $\{1,\dots, n+1\}$, let $\tau_i,\, i =
1,\dots , n$ be the transposition $(i,n+1)$ and let $\tau_{n+1}$
be the identity. Then every element of $S_{n+1}$ can be written
uniquely as a product $\sigma\tau_i$ for some $i$, where $\sigma$
is a permutation of $1,\dots, n$. Now let $C_2 = \{1,a\}$ and
view $C_2^n$ again as the subgroup of $C_2^{n+1}$ consisting of
all $(c_1,\dots,c_{n+1})$ where $c_i$ is either $1$ or $a$ and
the number of entries equal to $a$ is even. If $c_1,\dots,c_n$
are known then so is $c_{n+1}$, which we may therefore omit from
the notation. With this, we have $\tau_i(c_1,\dots,c_n) =
(c_1',\dots,c_n')$ where $c_j' = c_j$ for $j\ne i$ and $c_i' = 1$
if the number of $c_i$ equal to $a$ is even and $c_i' = a$ if it
is odd.   
Now $\Bbb QC_2 \cong \Bbb Q \oplus \Bbb Q$ where the idempotents
of the summands on the right side may be taken to be $e =
(1+a)/2$ and $f = (1-a)/2$, respectively. Using our previous
notation, the $2^n$ primitive idempotents of $\Bbb QC_2^n$ may
therefore be written in the form $E = e(I)$, where $I$ is an
ordered $n$-tuple of indices whose values can only be $1$ or $2$;
the $j$th entry in $e(I)$ is $e$ or $f$ according as the $j$th
entry in $I$ is $1$ or $2$. Since $S_{n+1}$ operates as
automorphisms of $\Bbb QC_2^n$ it must permute these idempotents
$e(I)$. The effect of a permutation $\sigma$ of $\{1,\dots,n\}$
is clearly just to permute the entries of $E$. We must compute 
$\tau_iE, i = 1,\dots,n$. For simplicity, suppose that $i=n$; the
argument will be the same for the other values of $i$. To avoid
confusion, write $u$ instead of $1$ for the unit of $C_2$, and
consider the non-commutative polynomial ring $\Bbb Q\{u, a\}$
formally generated by these two symbols. If we consider $e$ and
$f$ also as abstract symbols, this is the same ring as $\Bbb
Q\{e, f\}$ under the identifications $e = (u+a)/2, f = (u-a)/2$.
Now consider the effect of $\tau_n$ on elements of the form $\xi
e, \xi f$, where $\xi$ is a monomial in $a$ and $u$. If the
degree of $\xi$ in $a$ is even then $\tau_n(\xi a) = \xi a,
\tau(\xi u) = \xi u$, while if the degree of $\xi$ in $a$ is odd
then $\tau_n(\xi a) = \xi u, \tau_n(\xi u) = \xi a$. It follows
that whatever the degree of $\xi$ in $a$ we have
$\tau_n(\xi e) = \xi e$ while $\tau_n(\xi f) = \pm \xi f$
according as the degree of $\xi$ in $a$ is even or odd. 

It follows from the foregoing that if $\xi$ is an arbitrary
homogeneous polynomial of total degree $n-1$ in $a$ and $u$ then
$\tau_n(\xi e) = \xi e$. In particular, if $E'$ is a monomial of
degree $n-1$ in $e$ and $f$ then $\tau_n(E'e) = E'e$. In
$\tau_n(E'f)$, however, if $E'$ is expanded as a polynomial in
$a$ and $u$, then the signs of all the terms of odd degree in $a$
are reversed. This is the same thing as replacing every factor
$e$ in $E'$ by $f$ and every factor $f$ by $e$. Since the
analogous result clearly holds for all $i = 1,\dots,n$, we have
the following description of the action of $S_{n+1}$ on $\Bbb
QC_2^n$ when the latter is viewed as the set of homogeneous
polynomials of total degree $n$ in $\Bbb Q\{e,f\}$. (This is, of
course, essentially the same thing as the tensor algebra of $\Bbb
QC_2$, and the polynomials of degree $n$ are just $(\Bbb
QC_2)^{\otimes n}$.)

\proclaim{Theorem 11}Let $E \in (\Bbb QC_2)^n$ be a monomial of
degree $n$ in $e$ and $f$. If the $i$th factor in $E$ is $e$,
then $\tau_iE$ = $E$. If the $i$th factor in $E$ is $f$ then the
$i$th factor of $\tau_iE$ remains $f$ but in all other factors
$e$ and $f$ are interchanged.$\qed$ \endproclaim

In this theorem, the fact that the ground field is $\Bbb Q$ is of
no consequence; it was used only to show that we had constructed
a representation of $S_{n+1}$ as operators on the $n$th cartesian
power of a two-element set.  In purely combinatorial terms it can
be rephrased as follows.
\proclaim{Theorem 12} There is a permutation representation of
$S_{n+1}$ on the set of sequences of length $n$ of $0$'s and
$1$'s in which the subgroup $S_n$ of $S_{n+1}$ fixing $n+1$
operates by permutation of the entries, and the transpositions
$(i, n+1)$ are represented by operators $\tau_i$ as follows: If
the $i$th entry of a sequence $\frak s$ is $0$ then $\tau_i\frak
s = \frak s$. If the $i$th entry of $\frak s$ is $1$ then the
$i$th entry of $\tau_i\frak s$ remains $1$ but all other entries
are replaced by their complements (i.e., $0$'s are replaced by
$1$'s and $1$'s by $0$'s). $\qed$ \endproclaim
 
It is clear that the representation is faithful except in the
trivial case where $n=1$. The group consisting of all
permutations of the entries in the sequences together with all
operations which replace a single entry by its complement is the
full group of permutations of all $2^n$ of the sequences of $0$'s
and $1$'s.

We are now free to use any coefficient ring $k$. Referring back
to the Hecke algebra construction (\S4) we can take $k=\Bbb
Z_{q,2} = \Bbb Z[q, q^{-1}, 1/2_{q^2}]$ where $2_{q^2} =
1/(1+q^2)$ and eventually $q = 1+t$. Then $kC_2 \cong k\oplus k$
where the idempotents generating the summands are, respectively,
$e= (1+qa)/2_{q^2}, f = (q^2-qa)/2_{q^2}$. Note that when $q=1$
these are just our previous $e$ and $f$, and that as $2_{q^2} =
1+q^2 = 2+2t+t^2$ we can reduce modulo any rational prime, in
particular $p=2$, and the expressions remain meaningful. When
$q=1$ or equivalently, $t=0$, the operation of $S_{n+1}$ is just
the original operation on $\Bbb QC_2^n$, but the deformation of
$\Bbb ZC_2$ to its Hecke algebra and its induced deformation of
$\Bbb ZC_2^n$ is not equivariant with respect to the original
operation. It is so with respect to the subgroup $S_n$ which just
permutes the factors, which in effect is just the subgroup
leaving $n+1$ fixed, but not with respect to the $\tau_i$.
Nevertheless, the theorem effectively defines an operation of all
$S_{n+1}$, but that operation is a deformation of the original
one. 

It is important now to exhibit the deformation explicitly in
order to show what happens at $p=2$. To do so, we must write it
in terms of the basis $\{1, a\}$ for $kC_2$. Here $e$ and $f$ are
represented by the column vectors of the matrix $X = \frac
1{2_{q^2}}\pmatrix 1&q^2\\q&-q\endpmatrix$. Putting $1$ before
$a$  gives a lexicographic order to the $2^n$ basis elements of
$k(C_2^n) = (kC_2)^{\otimes n}$ consisting of the tensor products
of length $n$ of $1$'s and $a$'s. The tensor products of length
$n$ of the $e$'s and $f$'s arranged in lexicographic order are
then represented by the column vectors of the $2^n \times 2^n$
matrix $X^{\otimes n}$. To illustrate what happens it will be
sufficient to consider the case $n=2$, so we are tacitly dealing
with $D_3$. Write 
$$Y = X^{\otimes 2} = \frac 1{(2_{q^2})^2}\pmatrix 
1 & q^2 & q^2 & q^4\\ q & -q & q^3 & -q^3\\ q & q^3 & -q & -q^3\\
q^2 & -q^2 & -q^2 & q^2 \endpmatrix;
Y^{-1} = \pmatrix 1 & q & q & q^2\\ 1 & -
q^{-1} & q & -1\\ 1 & q & -q^{-1} &-1 \\ 1 & -q^{-1} & -q^{-1} &
q^{-2} \endpmatrix.
$$
If an element $\sigma \in S_3$ is represented in terms of the
basis $\{e\otimes e, e\otimes f, f\otimes e, f\otimes f\}$ by a
matrix $M$, then its representation in terms of of the basis
$\{1\otimes 1, 1\otimes a, a\otimes 1, a\otimes a\}$ is then
given by $YMY^{-1}$. For $\sigma = (1,2)$ the effect on this
basis is to interchange its second and third elements, so the
matrix $M$ is the permutation matrix 
$P_{23} = \pmatrix 1 & 0 & 0 & 0\\ 0 & 0 & 1 & 0 \\ 0 & 1 & 0 &
0\\ 0 & 0 & 0 & 1 \endpmatrix$. One can check immediately that
this commutes with $Y$, so the operation of $(1,2)$ is
undeformed. More generally, it is easy to see in the case of
$D_{n+1}$ that the operation of $S_n$, which operates by the
usual permutation of the tensor factors, is undeformed. The
transposition $(n, n+1)$, in our illustration $(2,3)$, behaves
quite differently. In the basis generated by the tensor products
of the idempotents $e,f$, it fixes $e\otimes e$ and $f\otimes e$
and interchanges $e \otimes f$ and $f\otimes e$, so it is
represented by $P_{24} = \pmatrix 1 & 0 & 0 & 0\\0 & 0 & 0 & 1\\
0 & 0 & 1 & 0\\ 0 & 1 & 0 & 0\endpmatrix$. Its matrix in terms of
the basis generated by the elements $1, a$ of $C_2$ is easily
calculated to be 
$$
YP_{24}Y^{-1} = \frac 1{(2_{q^2})^2} \pmatrix
(2_{q^2})^2 & 0 & q^5-q & 1-q^4\\
0 & (2_{q^2})^2 & 1-q^4 & q^3 - q^{-1}\\
0 & 0 & 1-q^4 & 2(q^3 + q) \\
0 & 0 & 2(q^3+q) & 1-q^4 \endpmatrix
$$
Letting $q \to 1$, or equivalently, writing $q = 1+t$ and letting
$t\to 0$, this becomes the permutation matrix $P_{34}$, which one
can check is indeed the representation of $(2,3)$ in the basis
generated by $1$ and $a$. This matrix remains formally unchanged
after reduction modulo any prime, in particular $p = 2$. However,
reduction modulo $2$ and reduction modulo $t$ do not commute
here. Reducing first modulo $2$ and then modulo $t$ yields the
matrix $N = \pmatrix 1 & 0 & 1 & 1\\ 0 & 1 & 1 & 1\\ 0 & 0 & 1 &
0 \\ 0 & 0 & 0 & 1 \endpmatrix$.   
Letting $q \to 1$ first recovers the original representation of
$(1,2)$ and $(2,3)$ on the $1,a$ basis by $P_{23}$ and $P_{34}$;
reducing first modulo $2$ gives, respectively, $P_{23}$ and the
above $N$. Setting $W = \pmatrix 1 & 1 & 1 & 1\\ 0 & 0 & 1 & 0\\0
& 1 & 0 & 0\\0 & 1 & 1 & 1 \endpmatrix$, one has $WP_{23}W^{-1} =
P_{23}, WNW^{-1} = P_{34}$, so $W$ conjugates the pair $P_{23},
N$ into $P_{23}, P_{24}$. Therefore, $P_{23}$ and $N$ give, a
representation of $S_3$ on $kC_2$ equivalent to that originally
given by $P_{23}$ and $P_{24}$. Thus, if we had started in characteristic $2$ we would
still have produced the necessary deformation of the action of
$S_3$ consistent with the deformation of $kC_2$ but its behavior
at different primes patches together in an unexpected way.
Nevertheless we do have a global deformation; it will behave
correctly at every prime.

Having now deformed the integral group ring $\Bbb ZD_n$ to a
smash product $A = (k\oplus k)^{\otimes n}\#kS_{n+1}$, Theorem 3
is applicable. In the next stage of the solution we must examine
the orbits and the isotropy groups of the action of $S_{n+1}$;
the decomposition into orbits will give a decomposition of the
algebra. If $B(I)$ is the direct sum of the rank one algebras
$A(J)$ for $J$ in the orbit of $I$, then we will again have $A =
\oplus B(I)\#k(S_{n+1})_I$ where the sum ranges over the orbits
of $S_{n+1}$ in the set of indices $I$. Since the $A(I)$ are all
just isomorphic to $k$ they have no automorphisms, so if the
orbit of $I$ has order $m$ then $B(I) \cong M_m(k)
\#k(S_{n+1})_I$.  For the Weyl groups $B_n = C_2\wr S_n$ the
isotropy groups were of the form $S_m \times S_{n-m}$. Since each
factor had a split global solution starting with $\Bbb Z$, so did
their product. Here the same is almost true, but for $D_{n+1}$
with $n$  odd there is a special case, the ``middle'' isotropy
group. This is the one case where we have to use Theorem 5.

To examine the orbits and isotropy groups for $D_{n+1}$ consider,
for simplicity, the case where $E = e(I)$ is a tensor product of
$m$ factors $e$ followed by $n-m$ factors $f$; 
write it as $e^mf^{n-m}$. Every $e(J)$ with the same number of
$e$'s and $f$'s as this is in the orbit of $E$. Generally this is
not the whole orbit since $\tau_nE$ will contain $e$ now 
$n-m-1$ times and $f$ now $m+1$ times. The exceptional or
middle case is that where $n = 2r+1$ is odd and $m=r$. 
One must be careful in counting orbits, since if we let $m$ run
from $0$ to $n$ then every orbit will be counted twice except the
middle one (when $n$ is odd). To count each only once, restrict
$m$ by requiring that $m \le n-m$ or $2m \le n$. The middle case
must also be handled separately.

In the non-exceptional case, the order of the orbit of $e(I)$ is
therefore $\binom nm + \binom n{m+1} = \binom {n+1}{m+1}$ and the
order of the isotropy group $(S_{n+1})_I$ is $(m+1)!(n-m)!$. In
fact, in the non-exceptional case, we have $(S_{n+1})_I \cong
S_{m+1} \times S_{n-m}$. Here the second factor arises because it
is evident that permuting the last $n-m$ factors of $e(I)$ will
leave it unchanged. For the first $m$ factors, since they are all
equal to $e$, we can not only permute them, but also apply
$\tau_i$ for $i=1,\dots, m$ and we have just proven that the
group generated by all these operations in $S_{m+1}$. (It may
appear that the case $m=n$ is also exceptional since $f$ does not
occur and the orbit of $e^m$ is reduced to $e^m$ itself, but the
formulas remain correct if we interpret $S_0$, whose order is
$0!$, as the identity.) In the exceptional case where $n = 2r+1$
and $e(I) = e^rf^{r+1}$, the order of the orbit is only $\binom
{2r+1}r$ and the order of the isotropy group is therefore
$(2r+2)!/\binom{2r+1}r = 2((r+1)!)^2$. This isotropy group
contains the subgroup isomorphic to $S_{r+1}\times S_{r+1}$,
where, as before, the first factor is generated by the
permutations of the first $r$ factors together with the $\tau_i$
for $i = 1,\dots,r$ and the second factor consists of the
permutations of the last $r+1$ factors. Being of index $2$, it is
normal.  

\proclaim{Lemma D} In the middle case, where $e(I) = e^rf^{r+1}$
one has $(S_{2r+2})_I \cong S_{r+1} \wr C_2$. \endproclaim
\demo{Proof} Set $\sigma = (1,r+1)(2,r+2)\cdots(r,2r)$ and $\rho
= \sigma\tau_{2r+1}$. Then $\rho e(I) = e(I)$ so $\rho$ is in the
isotropy group of $I$. Clearly $\rho^2 = 1$. What one must show
is that conjugation by $\rho$ interchanges the factors in the
product $S_{r+1} \times S_{r+1}$. It is trivial that $\rho$
interchanges the permutations of places $1,\dots,r$ with the
permutations of places $r+1,\dots,2r$. Finally, we claim that
$\rho\tau_i\rho = (r+i, 2r+1)$ for $i = 1,\dots,r$. One must show
that the operations are the same on every $e(J)$. There are four
cases, according as the $i$th entry of $e(J)$ is $e$ or $f$ and
likewise for the $(2r+1)$st entry; each is readily checked.
$\qed$
\enddemo
With this the Donald--Flanigan problem for $D_n$ is solved. For
with a split global solution for $S_{r+1}$, the problem for the
middle isotropy group is reduced to our previous one of a wreath
product of a sum of central separable algebras. The largest
isotropy group is $S_{n+1}$, so we have the following.

\proclaim{Theorem 13} There is a split global solution starting
with $\Bbb Z$ to the Donald--Flanigan problem for the Weyl groups
$D_{n+1}$; all their irreducible complex representations are
rational. The coefficient ring after deformation is $\Bbb
Z_{q,n+1}$. $\qed$
\endproclaim

We can also compute the rational group ring of $D_{n+1}$ in the
following way, separating the odd and even cases, and just
writing $n=2r$ in the latter. For the middle case we need to
compute $\Bbb Q(S_r\wr C_2)$. From Theorem 5, for each
block $S_r(\lambda)$ of $\Bbb QS_r$ there is a pair of summands
each isomorphic to $S_r(\lambda)$, and for every pair of distinct
blocks $S_r(\lambda), S_r(\mu)$ of $\Bbb QS_r$ there is a unique
summand isomorphic to $M_2(S_r(\lambda)\otimes S_r(\mu))$.  The
order of $\lambda$ and $\mu$ is immaterial; there is only one
summand for the pair. Write $\lambda < \mu$ to indicate that
$\lambda$ precedes $\mu$ in the lexicographic order of partitions
of $n$.
\proclaim{Lemma E} The rational group ring for the middle
component of $\Bbb QD_{2r}$ is given by
$$
\Bbb Q(S_r\wr C_2) =
\bigoplus_{\lambda}\{S_r(\lambda)\oplus S_r(\lambda)\}
\bigoplus\Sb \lambda,\mu \\ \lambda < \mu \endSb \{M_2(\Bbb Q)
\otimes S_r(\lambda)\otimes S_r(\mu)\} \qed
$$
\endproclaim
This gives the following expression for the rational group ring.
\proclaim{Theorem 14} The rational group ring of the Weyl groups
of type $D_n$ is given by
$$\align
\Bbb QD_{2r+1} &=\,\,
\bigoplus_{m=0}^{r}\,\{M_{\binom{2r+1}{m}}(\Bbb Q)\otimes \Bbb
QS_{2r+1-m} \otimes \Bbb QS_{m}\} \\
\Bbb QD_{2r}\quad &= \bigoplus_{m=0}^{r-1} \!\!
\{M_{\binom {2r}{m}}(\Bbb Q) \otimes \Bbb QS_{2r-m}\otimes \Bbb
QS_{m}\} \bigoplus\{M_{\binom {2r-1}r}(\Bbb Q) \otimes \Bbb Q(S_r
\wr C_2)\},
\endalign
$$
where $\Bbb Q(S_r\wr S_r)$ is given by Lemma E. $\qed$
\endproclaim
From this one can see that the representations of $D_n$ may be
indexed by unordered pairs of representations of $S_m$ and 
$S_{n-m}$. If $n$ is odd then the order of the pair will be fixed
by the requirement that $2m >n$ while if $n=2r$ is even then, as
we have just seen, it is true both for $m\ne r$ and for $m = r$.

As mentioned at the beginning, now only six finite reflection
groups remain for which a solution to the Donald--Flanigan
problem is unknown. If the problem is solvable for all of them,
one would like to know whether in all cases the Hecke algebra
effectively gives a solution. While the Donald--Flanigan problem
has been solved here for large classes of groups, the major
question remains of what characterizes a group for which it is
solvable. That, so far, is a mystery.

\enddocument